\input amstex
\documentstyle{amsppt}
\pageheight{23cm}
\pageheight{23cm}
\pagewidth{16.8cm}

\UseAMSsymbols 
\loadbold 
\centerline{\bf A Cohomology (p+1) Form Canonically Associated with Certain 
Codimension-q Foliations}\centerline{\bf on a Riemannian Manifold } 
\bigskip 
\bigskip 
\centerline{} 
\vskip .5in 
\centerline{\bf by Gabriel Baditoiu, Richard H. Escobales, Jr., and Stere Ianus} 
\vskip .5in 

\item{ }\underbar{Abstract:} Let $(M^{n},g)$ be a closed, connected, 
oriented, $C^{\infty}$, Riemannian, $n$-manifold with a transversely oriented foliation 
$\boldkey F$. We show that if $\lbrace X,Y \rbrace$ are basic vector fields, the leaf 
component of $[X,Y]$, $\Cal{V}[X,Y]$, has vanishing leaf divergence whenever 
 ${\kappa}\wedge \chi_{\boldkey F}$ is a closed (possibly zero) de Rham cohomology $(p+1)$-form. Here 
 ${\kappa}$ is the mean curvature one-form of the foliation ${\boldkey F}$ and 
 ${\chi_{\boldkey F}}$ is its characteristic form. ÊIn the codimension-$2$ case, 
 ${\kappa}\wedge \chi_{\boldkey F}$ Êis closed if and only if ${\kappa}$ 
 is horizontally closed. In certain restricted cases, we give necessary and sufficient conditions for 
 ${\kappa}\wedge{\chi_{\boldkey F}}$ to be harmonic. As an application, 
 we give a characterization of when certain closed $3$-manifolds
 are locally Riemannian products.
We show that bundle-like foliations with totally umbilical leaves with leaf 
dimension greater than or equal to two on a constant curvature manifold, with 
non-integrable transversal distribution, and with Einstein-like transversal 
geometry are totally geodesic.

\bigskip 
\noindent 
2000 {\it Mathematics  subject  classification (Amer. Math. Soc.)}: 
primary 57R30; 
secondary 53C25. 


\vskip .5in 
\noindent {\bf Introduction}
\vskip .5in 

\noindent 
Let $(M,g)$ be an oriented, $n$-dimensional Riemannian manifold admitting a 
transversely oriented
foliation, ${\boldkey F}$,  of leaf dimension $p$ and codimension $q$ so $p + 
q = n$. Generally, 
we will assume that both integers $p$ and $q$ are positive. Let ${\kappa}$ 
denote the mean 
curvature  one-form associated with the foliation ${\boldkey F}$ and let  
${\chi}_{\boldkey F}$ 
denote the characteristic form of the leaves  of ${\boldkey F}$. Following 
Kamber and Tondeur, 
we consider the  $p+1$ form ${\kappa} \wedge {\chi}_{\boldkey F}$. Suppose 
$X$ and $Y$ are local
basic  vector fields orthogonal to the leaves of ${\boldkey F}$. If 
${\kappa} \wedge {\chi}_{\boldkey F}$
is a closed form, then the leaf divergence of the leaf component of $[X,Y]$, 
$div_{\boldkey F} {\Cal V}[X,Y]$, vanishes identically. In the special case of 
$q=2$, 
${\kappa} \wedge {\chi}_{\boldkey F}$ is closed 
if and only if $div_{\boldkey F} {\Cal V}[X,Y] \equiv 0$ for any local basic vector 
fields $X$ and 
$Y$. This result, Theorem 1.2 below, is a general result for a foliation on an 
{\it arbitrary} 
Riemannian manifold. It illustrates once more the general principle that when a 
cohomology form arises, some pleasant geometric consequences often follow.  
\newline
\newline
Now suppose that $(M,g)$ above is a closed manifold and that the foliation 
${\boldkey F}$ is a 
Riemannian foliation. Then a fundamental result of Dominguez asserts that 
there then exists a metric 
$g$ on $M$ so that the associated mean curvature one form ${\kappa}$ of 
${\boldkey F}$ with respect
to $g$ is a  basic one-form. In this setting a result of Kamber-Tondeur asserts 
that ${\kappa}$
is a closed one-form. Suppose now $q=2$. Then it is easy to see that 
${\kappa} \wedge {\chi}_{\boldkey F}$ is a closed form. In fact, we establish in 
Theorem 1.4 that
 ${\kappa} \wedge {\chi}_{\boldkey F}$ is  co-closed if and only if $ 
div_{\boldkey H} {\tau} = 0 $,
where ${\tau}$ is the mean curvature  vector field dual to ${\kappa}$ and where 
the divergence is
taken with respect to a basic orthonormal frame orthogonal to ${\boldkey F}$. 
The proof involves 
lengthy and not entirely routine calculations, using three sets of arguments. The 
O'Neill tensors $T$
and $A$ play a crucial role. 
\newline
\newline
Applying Theorem 1.4 to a Riemannian flow on a closed $3$-manifold, $M^{3}$, 
we show that with
respect to a Dominguez metric, $M^3$ decomposes as a local Riemannian 
product if and only if
$Ric^M (V,V) = 0$, where $V$ is the unit length vector field tangent to the flow, 
and 
${\kappa} \wedge {\chi}_{\boldkey F}$ is harmonic (Corollary 1.5).  
The proof uses an important result of Ranjan as developed in [T3].  Using 
computations from [T3], we
show additionally in Corollary 1.6 that for a Riemannian flow on closed $M^{n}$ 
splits as a local 
Riemannian product with respect to the Dominguez metric, if and only if 
$Ric^M (V,V) = 0$ and 
$div_{\boldkey H}{\tau} = 0$. Since the result of Corollary 1.6 is general, the 
computation does
not depend on Theorem 1.4.
\newline
\newline
Theorem 1.7 establishes a result similar to that of Theorem 1.4 in the less 
interesting case
$q=1$, while Theorems 1.8, 1.9 and 1.11 wrap things up in the spirit of [T3].
\newline
\newline    
\noindent 
In section \S 2 we obtain local properties for bundle-like foliations 
with totally umbilical leaves on a constant curvature manifold.
In Proposition 2.3 we obtain an equivalent condition for 
$\kappa$ to be horizontally closed for a bundle-like foliation with totally 
umbilical leaves. Now we assume the $n$-dimensional  Riemannian manifold $(M,g)$ has constant 
curvature $c$. In Proposition 2.4 we show that $\kappa$ is a basic one-form.
As a consequence of Theorems 1.4 and 1.7, we get that 
${\kappa}\wedge\chi_{\boldkey F}$ is harmonic if and only if $g(\tau,\tau)=-pqc$
provided that the transversal distribution is integrable,
the dimension of the leaves $p$ is greater than one and the codimension 
$q$ is either one or two. Assuming that the transversal distribution is non-integrable
(at any point), we obtain a sufficient condition for $\boldkey F$ to be totally geodesic.   
In an important paper, Walschap showed that a bundle-like foliation with totally umbilical leaves 
and with leaf dimension 
$p>1$ on a complete simply connected space of constant curvature $c\geq 0$ is totally geodesic (see [Wa]). 
Using a different approach, Theorem 2.8 provides a similar result to Theorem 3.1 in [Wa], 
under no global assumptions and under some additional local ones. A remarkable fact is that 
Propositions 2.3, 2.4, 2.5, Theorem 2.8 can be extended to the pseudo-Riemannian case 
with definite induced metrics on leaves,
and, additionally only for Theorem 2.8, with induced positive definite transversal metrics.

\vskip .5in 
\noindent {\bf 1. The (p+1) form ${\kappa}\wedge{\chi}_{\boldkey F}$} 
\vskip .5in 
 
\noindent Throughout this paper all maps, functions and morphisms are assumed to 
be at least of class $C^{\infty}$. On a closed connected oriented 
$C^{\infty}$ Riemannian manifold Ê$(M^{n},g)$, let $\boldkey F$ be a 
transversely oriented foliation 
of leaf dimension $p$ and codimension $q = n-p$. 
Let $\boldkey V$ Êdenote the distribution 
tangent to the foliation $\boldkey F$, and $\boldkey H$ 
the distribution orthogonal to $\boldkey V$ Êin $TM$ determined by the 
metric $g$. If $E$ is a Êvector field on $M$, $\Cal{V}E$ and $\Cal{H}E$ 
will denote the projections of $E$ onto Êthe distributions $\boldkey V$ Ê 
and $\boldkey H$ Êrespectively. 
Call the vector field $E$ {\it vertical} Êif $\Cal{V}E = E$. Call $E$ {\it 
horizontal} if $\Cal{H}E = E$. 

\bigskip 
\noindent 
In general a $C^{\infty}$ {\it Êfoliation} of {\it codimension-$q$} on an 
$n$-dimensional manifold $M$ can be defined is a maximal Êfamily of 
$C^{\infty}$ submersions 
 $f_{\alpha}: U_{\alpha} \to f_{\alpha}(U_{\alpha}) \subset \Bbb 
R^{q}$ where $\{ U_{\alpha}\}_{\alpha Ê\in \Lambda }$ is an open cover of 
$M$ and where for each ${\alpha},{\beta} \in {\Lambda}$ and each 
$x \in U_{\alpha} \cap U_{\beta}$, there exists a local diffeomorphism 
${\phi}^{x}_{\beta\alpha}$ 
of $R^{q}$ so $f_{\beta}$ = Ê${\phi}^{x}_{\beta ,\alpha} \circ f_{\alpha}$ 
in some neighborhood $U_{x}$ of $x$ (see [L], 2,3). 

\medskip 
\noindent 
A horizontal vector field $Z$ defined on some open set $U$ where $U 
\subset U_{\alpha}$ is called {\it Ê$f_{\alpha}$- basic} Êprovided 
$f_{\alpha*}Z$ is a well defined vector Êfield on $f_{\alpha}(U)$. ÊAs 
pointed out in [E1] (for {\it any} 
metric $g$), if $U \subset U_{\beta}$, then $Z$ is also {\it 
$f_{\beta}$- basic}, so one can speak of $Z$ as a {\it local basic 
vector field}. ÊWe sometimes drop the word ``local.'' Ê 
Let $i(W) $ and $\theta(W)$ denote the interior product 
and the Lie derivative with respect to a vector field $W$. A 
differential form $\phi$ is called {\it basic} provided 
$ i(W) {\phi} = 0 $ 
and ${\theta(W)}{\phi} = 0$ for all vertical vector 
fields $W$ ([T1], 118). We follow the conventions of [AMR] 
for the formalism of differential forms and their exterior derivatives. 
\medskip 
\noindent 
$D$ will denote the Levi-Civita connection on $M$ and, following [EP], we 
introduce the tensors $T$ and $A$ as follows. For vector fields $E$ and 
$F$ on $M$, 

$$
\alignat 2 
&  &\qquad \qquad & T_{E}F = Ê\Cal{V}D_{\Cal{V}E} \Cal{H}F + Ê 
\Cal{H}D_{\Cal{V}E} \Cal{V}F, \text{ and}\tag{1.1}\\ 
 & & \qquad \qquad & A_{E}F = \Cal{V}D_{\Cal{H}E} \Cal{H}F + 
 \Cal{H}D_{\Cal{H}E} \Cal{V}F Ê.\tag{1.2}\\ 
\intertext{Then $T$ and $A$ are tensors of type $(1,2)$. These 
tensors satisfy the usual properties outlined in [EP]. 
 We note that if $X$ and $Y$ are 
horizontal,} 
 & & \qquad \qquad & A_{X}Y \neq Ê- A_{Y}X,\text{ in general,} \tag{1.3}\\ 
\endalignat $$ 
\noindent 
unless the foliation $\boldkey F$ Êis $bundle-like$ with respect 
to the metric $g$ (see [JW], Lemma (1.2)) that is, if $X$ is a basic 
vector field, 
$W g(X,X) = 0$ for every vertical vector field $W$. 
If $\lbrace V_{1}, V_{2}, V_{3},... V_{p} \rbrace$ is a local 
orthonormal frame tangent to the foliation, we define the {\it mean 
curvature one-form} Ê$\kappa$ as follows: 

$$ 
\leqno (1.4) \qquad \qquad \qquad\kappa(E) = \sum_{i = 1}^{p} g(E, T_{V_{i}}V_{i}). 
$$ 
\noindent 
Call $\kappa$ {\it horizontally closed} if $d \, \kappa(Z_{1}, 
Z_{2})=0$ for any horizontal vector fields $Z_{1}, Z_{2}$. Using the usual 
properties of the tensor $T$, one sees easily that if $X$ is basic, 

$$\leqno (1.5) Ê\qquad\qquad \qquad \kappa(X) = \sum_{i=1}^{p}g([X,V_{i}], V_{i}) .$$ 

\noindent 
Following [T-1], page 65--66, let $\chi_{\boldkey F}$ denote the 
characteristic form for the foliation $\boldkey F$. Then with $\lbrace V_{1},...,V_{p}\rbrace$ 
as above 
and for vector fields Ê$\lbrace E_{1},\dots, E_{p} \rbrace$ on $M^{n} $, we have: 

$$\leqno (1.6) Ê\qquad \qquad \chi_{\boldkey F}(E_{1}, E_{2,\ldots}, E_{p}) = 
det(g(E_{i},V_{j})) .$$ 
This characteristic differential form (see [T3], page 37) is independent of the 
local orthonormal frame Ê$\lbrace V_{1},... V_{p} \rbrace$. 
If any one of the arguments $E_{i}$ is horizontal, then the 
left hand side of (1.6) vanishes. This fact will be used repeatedly in 
the computations below. 
\newline 
\newline 
\noindent 
We say a ${\boldkey F}$ is a Ê{\it Riemannian} foliation of 
leaf dimension $p$ and codimension-$q$, provided that there is 
some Riemannian metric $g$ on $M^{n}$ with respect to which 
${\boldkey F}$ is bundle-like in the sense above. If 
${\boldkey F}$ is a Riemannian foliation on a compact manifold 
$M^{n}$, then a fundamental result of Dominguez, [D], shows 
that there Êalways exists a metric $g$ for which 
${\boldkey F}$ is bundle-like and for which the associated 
mean curvature one-form, ${\kappa}$, is basic. We call this 
metric, a {\it Dominguez metric}. 
\newline 
\newline 
\noindent 
Part (a) of the following result is proven in the appendix of [EP]. 
\bigskip 
\bigskip 
\noindent \underbar{Lemma 1.1.} Ê(a) Let $(M^{n},g)$ be a connected, 
oriented, $C^{\infty}$ ÊRiemannian $n$-manifold with a transversely oriented 
codimension-$q$ foliation $\boldkey F$, with $q \ge 2$. 
Suppose $X$ and $Y$ are basic vector fields. Then $\Cal{V}[X,Y]$ has vanishing leaf divergence 
if and only if ${\kappa}$ is horizontally closed. 
\newline\noindent 
(b) Let ${\boldkey F}$ be a transversely oriented 
Riemannian foliation on a closed, oriented Riemannian manifold $(M,g')$. 
Then there exists a Dominguez metric $g$ on $M$ so that if $X$ and $Y$ are 
basic with respect to $g$, then $div_{\boldkey F} {\Cal V}[X, Y] = 0$ and 
indeed $div_{M}{\Cal V}[X,Y] = 0$. 
\newline 
\newline 
\noindent \underbar{Proof}. (a) This follows immediately from formula (3) of [EP] which can be expressed 
this way: 

$$\leqno (1.7) Ê\qquad \qquad d\kappa(X,Y) = - div_{\boldkey F} {\Cal V}[X,Y] ,$$ 

\noindent where the right hand side denotes the divergence of ${\Cal V}[X,Y]$ along a leaf of ${\boldkey F}$. 
A more succinct proof of (1.7) Êappears in [CE 1]. 
\medskip 
\noindent 
(b) Let $g$ be such a metric for $M$. Then ${\boldkey F}$ 
is bundle-like with respect to $g$ and the associated mean curvature one 
form ${\kappa}$ is basic by Dominguez's Theorem. Then ${\kappa}$ is 
closed by a result of Kamber-Tondeur [T3, p. 82], and so in particular, ${\kappa}$ 
is horizontally closed. Thus, by Lemma 1.1(a) Êand the appendix to [EP], 
$div_{\boldkey F} {\Cal V}[X, Y] = div_{M}{\Cal V}[X,Y] = 0 $.

\bigskip 
\bigskip 
\noindent 
The form ${\kappa}\wedge \chi_{\boldkey F}$ arises in the important role in the work of Kamber 
and Tondeur on foliations, especially Riemannian foliations ([T1], pages 121 and 152, [T3], page 82). ÊIt turns out that when
this form is closed, the following pleasant property obtains for {\it arbitrary} foliations on Riemannian manifolds of 
codimension $q \ge 2$ (actually $q \ge 1$). The result illustrates once more the tie between cohomology and geometry.
\bigskip 
\bigskip 
\noindent \underbar{Theorem 1.2.} Let $(M^{n},g)$ be a closed, connected, 
oriented, 
$C^{\infty}$ Riemannian $n$-manifold with a transversely oriented 
codimension-q foliation $\boldkey F$. 
Suppose $X$ and $Y$ are basic vector fields. Then $\Cal{V}[X,Y]$ has vanishing leaf divergence 
(equivalently ${\kappa}$ is horizontally closed) whenever 
${\kappa}\wedge \chi_{\boldkey F}$ is a closed (possibly zero) de Rham 
cohomology $p+1$ form. In fact, if the codimension of ${\boldkey F}$, $q=2$, then 
${\kappa}$ is horizontally closed if and only if ${\kappa}\wedge \chi_{\boldkey F}$ is closed. 

\bigskip 
\bigskip 

\noindent 
\underbar{Proof.} Note, if $q=1$, ${\kappa}\wedge{\chi_{\boldkey F}}$ 
is an $n$-form and hence closed. If $X$ and $Y$ are basic vector fields, 
then in the codimension-one case, $X = fZ$ and $Y = hZ$ where $f$ and $h$ are functions defined 
on an appropriate open set and $Z$ is a unit length horizontal vector field on that set. 
Then, ${\Cal V}[X,Y] = 0$ and the theorem always holds in this trivial case. 
\newline 
\newline 
\noindent 
To establish this for $q \ge 2$, we will use the local frame, 
$\lbrace V_{1}, V_{2}, V_{3},... V_{p}, X,Y \rbrace$ 
where $X$ and $Y$ are basic and span Ê$\boldkey H$ at each 
$x\in U$ where $X_{x}$ and $Y_{x}$ are defined. Note, 
we make no requirements that $X$ and $Y$ form a basic {\it orthonormal} 
frame for $\boldkey H$, since we do not yet assume the metric $g$ on $M$ is bundle-like. 
A fundamental result of Rummler [Ru], yields: 

$$\qquad \qquad d{\chi_{\boldkey F}}(V_{1},\ldots ,V_{p}, X) = 
(-1)^{p+1}{\kappa}(X){\chi_{\boldkey F}}(V_{1},\ldots, V_{p}) . \tag{1.8}$$ 
\newline 
\newline 
\noindent 
At this point it is worth pointing out that for any $(p+2)$-form ${\gamma}$, 
$${\gamma}(V_{1}, V_{2}, V_{3},... V_{p}, X,Y) = {\gamma}(X,Y,V_{1}, V_{2}, 
V_{3},... V_{p}).$$ 
Using the formulas in [AMR] page 394 and the 
remarks above, we have, 

$$
\alignat{3}
& d({\kappa}\wedge \chi_{\boldkey F})(V_{1},\ldots,V_{p}, X,Y) &=& 
(d{\kappa}\wedge \chi_{\boldkey F})(V_{1},\ldots,V_{p},X,Y) - ({\kappa}\wedge 
d\chi_{\boldkey F})(V_{1},\ldots,V_{p}, X,Y)  \tag{1.9}\\
&  &=& (d{\kappa}\wedge \chi_{\boldkey F})(X,Y,V_{1},\ldots,V_{p}) 
- ({\kappa}\wedge d\chi_{\boldkey F})(V_{1},\ldots,V_{p},X,Y),
\endalignat
$$ 
\newline 
\noindent 
 which becomes, 
\newline 
\noindent 
$$\alignat{2} (1.10)
\qquad
&=&d\kappa(X,Y) \chi_{\boldkey 
F}(V_{1},\ldots,V_{p}) - (-1)^{p}\kappa(X)d\chi_{\boldkey F}(V_{1},\ldots,V_{p},Y) - 
(-1)^{p+1}\kappa(Y)d\chi_{\boldkey F}(V_{1},\ldots,V_{p},X)\quad
\\
&=&d\kappa(X,Y)\chi_{\boldkey F}(V_{1},\ldots,V_{p}) + (-1)^{p+1}\kappa(X)d\chi_{\boldkey 
F}(V_{1},\ldots,V_{p},Y) + (-1)^{p+2}\kappa(Y)d\chi_{\boldkey 
F}(V_{1},\ldots,V_{p}, X),
\endalignat
$$ 
\newline 
\noindent 
which then becomes by (1.8), 
\newline 
\newline 
\noindent 
$$\leqno (1.11) \qquad d\kappa(X,Y)\chi_{\boldkey F}(V_{1},\ldots,V_{p}) 
+ (-1)^{p+1}\kappa(X)(-1)^{p+1}\kappa(Y)\chi_{\boldkey 
F}(V_{1},\ldots,V_{p})$$ 
$$ \leqno \qquad\qquad\quad + (-1)^{p+2}\kappa(Y)(-1)^{p+1}\kappa(X)\chi_{\boldkey 
F}(V_{1},\ldots,V_{p})$$ 
\newline 
\newline 
\noindent which becomes, 
$$ \leqno (1.12)\qquad\qquad d\kappa(X,Y), $$ 
\newline 
\newline 
\noindent 
since $\chi_{\boldkey F}(V_{1},\ldots,V_{p}) = 1 $. Thus, 

$$ \leqno (1.13)\, Ê\qquad \qquad ÊÊd({\kappa}\wedge 
\chi_{\boldkey F})(V_{1},\ldots,V_{p},X,Y) = Êd\kappa(X,Y).$$ 
\newline 
\newline 
\noindent Hence, if ${\kappa}\wedge{\chi_{\boldkey F}}$ is a closed $(p+1)$-form, then 
in particular, the left hand side of $(1.13)$ vanishes, so ${\kappa}$ is 
horizontally closed and the result follows by Lemma 1.1. 
In the codimension-$2$ case, 
${\kappa}\wedge{\chi_{\boldkey F}}$ is closed if and only if 
the left hand side of $(1.13)$ vanishes. 
The proof of Theorem 1.2 is now complete.
\newline 
\newline 
\noindent 
We offer the following improvement of a result that appeared in [E2]. It should be noted that in Theorem 1.3 below 
we do {\it not} require that the flow ${\boldkey F}$ is {\it Riemannian} in the sense above. 
\bigskip 
\noindent 
\underbar{Theorem 1.3.} Let $(M^{3},g)$ be a closed, connected, oriented $C^{\infty}$ 
Riemannian manifold of dimension $3$ with a transversely oriented flow ${\boldkey F}$. 
Suppose the following conditions obtain. 
\newline 
\noindent 
(a) ${\boldkey F}$ admits a basic transverse volume form ${\mu}$. 
\newline 
\noindent 
(b) Ê${\kappa}\wedge{\chi_{\boldkey F}}$ is a closed $2$-form. 
\newline 
\noindent 
(c) Let $X$ and $Y$ denote local basic vector fields so that ${\mu}(X,Y) = 1$. 
Assume the globally defined vector field, ${\Cal V}[X,Y]$ satisfies the following: 

$[Z, {\Cal V}[X,Y]] = f_{Z} {\Cal V}[X,Y]$ 
\newline 
\noindent 
for any basic vector field $Z$ and for some function $f_{Z}$ depending on $Z$. 
\newline 
\noindent 
Then either: 
\newline 
\noindent 
 i) ${\Cal V}[X,Y]$ vanishes identically on $M$, so $\boldkey H$ is integrable 
and the leaves of ${\boldkey H}$ are minimal surfaces in $M$, or, 
\newline 
\noindent 
ii) ${\Cal V}[X,Y]$ never vanishes and so ${\boldkey H}$ is always a contact structure. 
\bigskip 
\noindent 
\underbar{Proof.} Condition (b) replaces the condition in [E2] that ${\kappa}$ is horizontally closed. 
These are equivalent, since when $n=3$, $p = 1$ and so $q=2$. 
In this case the last part of Theorem 1.2, (in particular (1.13)), 
applies so ${\kappa}$ is horizontally closed if and only if 
${\kappa}\wedge{\chi_{\boldkey F}}$ Êis closed. Then 
the argument given in [E2] carries over and ${\boldkey H}$ is a foliation of M by minimal surfaces, whenever 
${\Cal V}[X,Y]$ vanishes at one point. The only other possibility is that ${\Cal V}[X,Y]$ never vanishes, 
and in this case, ${\boldkey H}$ is a contact structure. 
\bigskip 
\bigskip 
\noindent 
Now assume that ${\boldkey F}$ is a {\it Riemannian} foliation of leaf dimension p and codimension $q = 2$. 
Although far more 
restrictive, this assumption allows us to give a definitive answer to the question: when is Ê 
${\kappa}\wedge {\chi}_{\boldkey F}$ co-closed? To address this problem, 
we need the following preparation. Because ${\boldkey F}$ is 
bundle-like with respect to $g$, the local submersions defining 
${\boldkey F}$ are Riemannian submersions in the sense of [O'N] (see 
[T1], [T2], [E] and [EP]), and so Êwe can choose a local orthonormal 
frame. Ê 
$\lbrace X_{1}, X_{2}, V_{1},V_{2},...,V_{p} \rbrace $ so $X_{1}, X_{2}$ are basic vector fields, and Êso 
$\lbrace V_{1}, V_{2},...,V_{p} \rbrace$ is a local orthonormal frame for ${\boldkey V}$. Indeed, 
at a fixed $x \in M$, we can choose Ê$X_{1}, X_{2}$ so Ê$ (\Cal{H} D_{X_{i}} X_{j})_{x} = 0$ for $ 1 \le i,j \le 2 $. 
Because, ${\boldkey F}$is Êbundle-like with respect to $g$, Ê$A_{X_{i}}X_{i} = 0$, for $1 \le i \le 2$. Ê 
Let $\lbrace V_{1}, V_{2},...,V_{p}\rbrace$ be a local orthonormal frame for ${\boldkey V}$. At a given $x \in M$, 
we can Êchoose this frame 
so that $(\Cal{V}D_{V_{i}}V_{j})_{x} = 0$. Note, for {\it any} vector Êfield $E$, 
$Eg(V_{i}, V_{i}) = 2g(D_{E}V_{i},V_{i})=0$ and $0 = Eg(V_{i}, V_{j}) = g(D_{E}V_{i}, V_{j}) + g(V_{i}, D_{E}V_{j})$. 
We will exploit these well known facts extensively in the computations below. When the orthonormal frame 
$\lbrace X_{1}, X_{2}, V_{1},V_{2},...,V_{p}\rbrace$ enjoys these additional properties at $x\in M$, we'll call the frame 
{\it a preferred orthonormal frame at x}. 
\newline 
\newline 
\noindent 
Set $ {\tau} = Ê\sum_{i = 1}^{p} \Cal{H} D_{V_{i}} V_{i} $. 
(We follow the conventions in foliations and suppress the usual constant.) Note, 
when $g$ is a Dominguez metric for ${\boldkey F}$, ${\tau}$ is basic, 
because it is dual to the one-form ${\kappa}$ with respect to $g$ by (1.4).
 Set $ div_{\boldkey H}{\tau} = g(X_{1}, D_{X_{1}}{\tau}) + g(X_{2}, D_{X_{2}} 
{\tau})$. ÊThen, Êit is well known ([P], page 151) that Ê$ \delta 
({\kappa}\wedge {\chi}_{\boldkey F}) = - div({\kappa}\wedge{ \chi}_{\boldkey F})$. 
We have the following result. 
Note, the operator ${\delta}$ below is on $M$ itself. 
\newline 
\newline 
\noindent 
\underbar{Theorem 1.4}. ÊLet $M^{n}$ be a closed, oriented, $C^{\infty}$, 
Riemannian manifold with a transversely oriented codimen\-sion-$2$ Riemannian foliation Ê 
${\boldkey F}$. Let $g$ be a Dominguez metric (with respect to which 
${\boldkey F}$ is bundle-like and ${\kappa}$ is basic). Then, ${\kappa}\wedge {\chi}_{\boldkey F}$ 
is closed. In fact, ${\kappa}\wedge{ \chi}_{\boldkey F}$ is harmonic if and only if 
$ div_{\boldkey H} {\tau} = 0$, where ${\tau}$ is the mean curvature vector field dual to ${\kappa}$. 
\newline 
\newline 
\underbar {Proof}. To show ${\kappa}\wedge{ \chi}_{\boldkey F}$ is closed, it suffices to observe 
that in the case under consideration, $d({\kappa}\wedge{ \chi}_{\boldkey F})$ is a closed $n$-form. 
Evaluating, $d({\kappa}\wedge{ \chi}_{\boldkey F})$ on a preferred orthonormal frame 
$ \lbrace X_{1}, X_{2}, V_{1},V_{2},...,V_{p} \rbrace$, we see that $(1.13)$ becomes 

$$\leqno (1.14) \qquad \qquad ÊÊd({\kappa}\wedge 
\chi_{\boldkey F})(V_{1},\ldots,V_{p},X_{1},X_{2}) = Êd\kappa(X_{1},X_{2}) .$$ 
\newline 
\noindent 
But, since ${\kappa}$ is basic, it is closed by a fundamental result of Kamber-Tondeur for bundle-like foliations 
with ${\kappa}$ basic (see [T3], page 82). This means that the left hand side of (1.14) vanishes and so 
${\kappa}\wedge{ \chi}_{\boldkey F}$ is closed, as claimed. 
\newline 
\newline 
\noindent We now will show that under the stated hypotheses, ${\kappa}\wedge{ \chi}_{\boldkey F}$ 
is co-closed on $M$. To do this we will use exclusively the preferred orthonormal frame above. 
The result follows from lengthy computations of $\delta ({\kappa}\wedge{ \chi}_{\boldkey F})$ 
on three sets of arguments: $(V_{1},V_{2},...,V_{p})$, $(X_{1},V_{1},\dots,\hat{V_{j}},\dots,V_{p})$ 
and $(X_{1},X_{2},V_{1}...,V_{p-2})$. It should be mentioned (that up to sign) it suffices to use $X_{1}$ 
in the second set of arguments. 
\newline 
\newline 
\noindent 
In the first evaluation, the reader should keep in mind the following
principles: 
$ g({\Cal V}D_{X_{i}}V_{j},V_{j})=0$, 
mentioned before. Secondly, terms with {\it repeated} vertical vector fields vanish, thirdly, the sum 
\newline 
\noindent 
$\sum_{i=1}^{p}({\kappa}\wedge{\chi}_{\boldkey{F}})(D_{V_{i}}V_{i},V_{1},V_{2},
...,V_{p})$ 
and the terms 
$({\kappa}\wedge{\chi}_{\boldkey F})(V_{1}, D_{V_{1}}V_{1}, V_{2},...,V_{p}),$ 
\newline 
\noindent 
$({\kappa}\wedge{\chi_{\boldkey F}})(V_{2},V_{1},D_{V_{2}}V_{2},...,V_{p}),\dots,$ 
 $({\kappa}\wedge{\chi_{\boldkey F}})(V_{p},V_{1},...,D_{V_{p}}V_{p})$ sum to zero. Also, 
$({\kappa}\wedge{ \chi}_{\boldkey F})$ vanishes identically on $(E_{1}, E_{2}, E_{3},\dots,E_{p+1})$ if 
two or more of the arguments are horizontal vector fields. 
Finally, Ê${\chi_{\boldkey F}}(E_{1},\dots,E_{p})$ will vanish identically 
even if all the $E_{j}$ are vertical but linearly dependent. In the expansion below, we have rearranged 
some of the terms in the expansion, but they are all there. We carry out each of these computations at $x\in M$ 
above. 
\bigskip 
$$ 
\alignat{3}
 &\delta (\kappa \wedge \chi_{\boldkey F})(V_{1},V_{2},\dots,V_{p}) &=& 
- \sum_{i = 1}^{p}(D_{V_i}(\kappa \wedge \chi_{\boldkey F}))(V_{i},V_{1},V_{2},\dots,V_{p})  \tag{1.15} 
\\
&&& -\sum_{a = 1}^{2} (D_{X_{a}} (\kappa \wedge \chi_{\boldkey F}))(X_{a},V_{1},V_{2},\dots,V_{p})
\endalignat 
$$
$$
\alignat{1}
=& - \sum_{i = 1}^{p}V_{i} ({\kappa}\wedge{\chi}_{\boldkey F})(V_{i},V_{1},V_{2},\dots,V_{p})
+   \sum_{i=1}^{p}({\kappa}\wedge{\chi}_{\boldkey F})(D_{V_{i}}V_{i},V_{1},V_{2},\dots,V_{p})\\ 
&+   \sum_{i=1}^{p}({\kappa}\wedge{\chi}_{\boldkey F})(V_{i}, D_{V_{i}}V_{1}, V_{2},\dots,V_{p}) 
+   \sum_{i=1}^{p}({\kappa}\wedge{\chi_{\boldkey F}})(V_{i}, V_{1},D_{V_{i}}V_{2},\dots,V_{p}) \\
&+ \dots + \sum_{i=1}^{p} ({\kappa}\wedge{\chi_{\boldkey F}})(V_{i}, V_{1},\dots,D_{V_{i}}V_{p})\\
%
%
%
%
%
&- \sum_{a = 1}^{2}X_{a}({\kappa}\wedge{\chi_{\boldkey F}})(X_{a},V_{1},V_{2},\dots,V_{p}) 
+ \sum_{a = 1}^{2} ({\kappa}\wedge{\chi_{\boldkey F}})(D_{X_{a}}X_{a},V_{1},V_{2},\dots,V_{p}) \\
&+ \sum_{a = 1}^{2}({\kappa}\wedge{\chi_{\boldkey F}})(X_{a},D_{X_{a}}V_{1},V_{2},\dots,V_{p})
+ \sum_{a = 1}^{2}({\kappa}\wedge{\chi_{\boldkey F}})(X_{a},V_{1},D_{X_{a}}V_{2},\dots,V_{p})\\
&+ \dots +\sum_{a = 1}^{2}  ({\kappa}\wedge{\chi_{\boldkey F}})(X_{a},V_{1},V_{2},\dots,D_{X_{a}}V_{p}).\\
\endalignat
$$ 
%
%
%
%
%
%
%
\newline 
\newline 
\newline 
\noindent Now Ê$(\Cal{H} D_{X_{i}}X_{i})_{x} = 0$, and 
$\Cal{V} D_{X_{i}}X_{i} = 0$ where defined, since the metric $g$ is assumed bundle-like. 
Since $\kappa$ annihilates vertical vector fields, $\chi_{\boldkey F}$ 
annihilates horizontal fields, the above becomes, 

$$
\alignat{3}
& & =& - X_{1} \kappa (X_{1}) - X_{2} \kappa 
(X_{2}) = - X_{1} g(X_{1}, \tau) - X_{2}g(X_{2}, \tau)
\tag{1.16}\\
& &=& - g(X_{1}, D_{X_{1}} \tau)-
g(X_{2},D_{X_{2}}\tau) = - div_{\boldkey H}\tau .
\endalignat 
$$ 
\newline 
\newline 
\noindent 
In the next expansion, note that $X_{1}$ is basic and that $\hat V_{j}$ means that $V_{j}$ is omitted. We can use $X_{1}$ 
as our basic vector field essentially without loss of generality. 

$$ \ \delta(\kappa \wedge 
\chi_{\boldkey F})(X_{1},V_{1},V_{2},\dots,\hat V_{j},\dots,V_{p}) \tag{1.17}$$ 
$$= -\sum_{a=1}^2 (D_{X_{a}}(\kappa \wedge \chi_{\boldkey F}))(X_{a},X_{1},V_{1},V_{2},\dots,\hat V_{j},\dots,V_{p}))$$ Ê 
 $$ - \sum_{i=1, i\neq j}^{p}(D_{V_{i}}(\kappa \wedge 
\chi_{\boldkey F}))(V_{i},X_{1},V_{1},V_{1},V_{2},\dots,\hat V_{j},\dots,V_{p}))$$ 
$$-(D_{V_{j}}({\kappa}\wedge{\chi}_{\boldkey F})(V_{j},X_{1},V_{1}, V_{2},\dots,\hat V_{j},\dots,V_{p})). $$ 
\bigskip 
\noindent 
Expanding (1.17), we have the following expression. 

$$- \sum_{a=1}^2 X_{a} (\kappa \wedge \chi_{\boldkey F})(X_{a},X_{1},V_{1},V_{2},\dots,\hat V_{j},\dots,V_{p}) \tag{1.18}$$
$$ + \sum_{a=1}^2  (\kappa \wedge \chi_{\boldkey F})(D_{X_{a}}X_{a},X_{1},V_{1},V_{2},\dots,\hat V_{j},\dots,V_{p}) $$ 
$$ + \sum_{a=1}^2 (\kappa \wedge \chi_{\boldkey F})(X_{a},D_{X_{a}}X_{1},V_{1},V_{2},\dots,\hat V_{j},\dots,V_{p}) $$ 
$$+ \sum_{a=1}^2  (\kappa \wedge \chi_{\boldkey F})(X_{a},X_{1},D_{X_{a}}V_{1},V_{2},\dots,\hat V_{j},\dots,V_{p}) +\dots $$ 
$$+ \sum_{a=1}^2 (\kappa \wedge \chi_{\boldkey F})(X_{a},X_{1},V_{1},V_{2},\dots,\hat V_{j},\dots,D_{X_{a}}V_{p}) $$ 
$$- \sum_{i=1, i\neq j}^{p}V_{i}({\kappa}\wedge{\chi_{\boldkey F}})(V_{i}, X_{1},V_{1},V_{2},\dots,\hat V_{j},\dots,V_{p})$$ 
$$ + \sum_{i=1, i\neq j}^{p}({\kappa}\wedge{\chi_{\boldkey F}})(D_{V_{i}}V_{i}, X_{1},V_{1},V_{2},\dots,\hat V_{j},\dots,V_{p})$$ 
$$ + \sum_{i=1, i\neq j}^{p}({\kappa}\wedge{\chi_{\boldkey F}})(V_{i},D_{V_{i}}X_{1},V_{1},V_{2},\dots,\hat V_{j},\dots,V_{p})$$ 
$$+  \sum_{i=1, i\neq j}^{p} ({\kappa}\wedge{\chi_{\boldkey F}})(V_{i},X_{1},D_{V_{i}}V_{1},V_{2},\dots,\hat V_{j},\dots,V_{p})$$ 
$$+  \sum_{i=1, i\neq j}^{p}({\kappa}\wedge{\chi_{\boldkey F}})(V_{i},X_{1},V_{1},D_{V_{i}}V_{2},\dots,\hat V_{j},\dots,V_{p})$$ 
$$+  \dots 
+ \sum_{i=1, i\neq j}^{p}({\kappa}\wedge{\chi_{\boldkey F}})(V_{i},X_{1},V_{1},V_{2},\dots,\hat V_{j},\dots,D_{V_{i}}V_{p})$$ 
%
%
%
%
$$- V_{j}({\kappa}\wedge{\chi_{\boldkey F}})(V_{j},X_{1},V_{1},V_{2},\dots,\hat V_{j},\dots,V_{p})$$ 
$$+ ({\kappa}\wedge{\chi_{\boldkey F}})(D_{V_{j}}V_{j},X_{1},V_{1},V_{2},\dots,\hat V_{j},\dots,V_{p})$$ 
$$+ ({\kappa}\wedge{\chi_{\boldkey F}})(V_{j},D_{V_{j}}X_{1},V_{1},V_{2},\dots,\hat V_{j},\dots,V_{p})$$ 
$$+ ({\kappa}\wedge{\chi_{\boldkey F}})(V_{j},X_{1},D_{V_{j}}V_{1},V_{2},\dots,\hat V_{j},\dots,V_{p})$$ 
$$+ ({\kappa}\wedge{\chi_{\boldkey F}})(V_{j},X_{1},V_{1},D_{V_{j}}V_{2},\dots,\hat V_{j},\dots,V_{p})$$ 
$$+ \dots + ({\kappa}\wedge{\chi_{\boldkey F}})(V_{j},X_{1},V_{1},V_{2},\dots,\hat V_{j},\dots,D_{V_{j}}V_{p}).$$ 

\bigskip 
\noindent
Most of the terms in (1.18) vanish for one of the following reasons: 
two of the arguments are horizontal; two repeated arguments. Note, at x, $D_{V_{i}}V_{k}$ is purely horizontal. 
The only non-zero summands in (1.18) are:
\bigskip
$$  (\kappa \wedge \chi_{\boldkey F})(X_{2},D_{X_{2}}X_{1},V_{1},V_{2},\dots,\hat V_{j},\dots,V_{p}) $$ 
$$- V_{j}({\kappa}\wedge{\chi_{\boldkey F}})(V_{j},X_{1},V_{1},V_{2},\dots,\hat V_{j},\dots,V_{p}) $$ 
$$+ ({\kappa}\wedge{\chi_{\boldkey F}})(V_{j},D_{V_{j}}X_{1},V_{1},V_{2},\dots,\hat V_{j},\dots,V_{p}).$$

\bigskip 
\noindent Since ${\Cal H}D_{V_{j}}X_{1} = A_{X_{1}}V_{j}$, this becomes: 
$$ Ê{\kappa}(X_{2}){\chi_{\boldkey F}}(A_{X_{2}}X_{1}, V_{1},V_{2},\dots,\hat{V_{j}},\dots,V_{p})\tag{1.19}$$ 
$$+ V_{j}({\kappa}\wedge{\chi_{\boldkey F}})(X_{1},V_{j},V_{1},V_{2},\dots,\hat{V_{j}},\dots,V_{p})$$ 
$$- ({\kappa}\wedge{\chi_{\boldkey F}})(A_{X_{1}}V_{j},V_{j},V_{1},V_{2}, \dots,\hat{V_{j}},\dots,V_{p}).$$ 
\bigskip 
\noindent 
Recall, $g(A_{X_{2}}X_{1},V_{j})V_{j} = g(X_{2}, A_{X_{1}}V_{j})V_{j}$. Now ${\tau}$ is basic, because ${\kappa}$ 
is basic and $g$ is bundle-like. Hence $ {\tau} = a_{1}X_{1} + a_{2}X_{2}$. Then (1.19) becomes: 

$$\ Ê{\kappa}(X_{2})g(A_{X_{1}}V_{j},X_{2}){\chi_{\boldkey F}}(V_{j},V_{1},V_{2},\dots,
\hat V_{j},\dots,V_{p}) \tag{1.20}$$ 
$$\pm V_{j}{\kappa}(X_{1}){\chi}(V_{1},V_{2},\dots,V_{j},\dots,V_{p})$$ 
$$-{\kappa}(A_{X_{1}}V_{j}){\chi_{\boldkey F}}(V_{j},V_{1},V_{2},\dots,\hat V_{j},\dots V_{p})$$ 
\bigskip 
\noindent which becomes 
$$ a_{2}g(A_{X_{1}}V_{j},X_{2}){\chi_{\boldkey F}}(V_{j},V_{1},V_{2},\dots,\hat V_{j},\dots,V_{p})
\tag{1.21}$$ 
$$\pm d{\kappa}(V_{j},X_{1}) - a_{2} g(A_{X_{1}}V_{j},X_{2})
{\chi_{\boldkey F}}(V_{j},V_{1},V_{2},\dots,\hat V_{j},\dots,V_{p}) = 0, $$ 
\newline 
because ${\kappa}$ is closed for the Dominguez metric. 
\bigskip 
\noindent 
Our final computation will involve evaluating $\delta ({\kappa}\wedge{\chi_{\boldkey F}})$ on 
$(X_{1},X_{2},V_{1},\dots,V_{p-2})$. Again, we can make this evaluation on our two basic fields and excluding 
$V_{p-1}$ and $V_{p}$ as arguments, essentially without loss of generality. 

$$ \delta({\kappa}\wedge{\chi_{\boldkey F}}))(X_{1},X_{2},V_{1},V_{2},\dots,V_{p-2})\tag{1.22}$$ 
$$= -\sum_{i=1}^{2}(D_{X_{i}}({\kappa}\wedge{\chi_{\boldkey F}}))(X_{i},X_{1},X_{2},V_{1},\dots,V_{p-2})$$ 
$$-\sum_{i=1}^{p}(D_{V_{i}}(({\kappa}\wedge{\chi_{\boldkey F}}))(V_{i},X_{1},X_{2},V_{1},V_{2},\dots,V_{p-2}).$$ 
\newline 
\noindent 
This expands to: 
$$ - \sum_{a=1}^{2}X_{a}({\kappa}\wedge{\chi_{\boldkey F}})(X_{a},X_{1},X_{2},V_{1},\dots,V_{p-2}) \tag{1.23}$$ 
$$ + \sum_{a=1}^{2}({\kappa}\wedge{\chi_{\boldkey F}})(D_{X_{a}}X_{a},X_{1},X_{2},V_{1},\dots,V_{p-2})
+ \sum_{a=1}^{2}({\kappa}\wedge{\chi_{\boldkey F}})(X_{a},D_{X_{a}}X_{1},X_{2},V_{1},\dots,V_{p-2})$$ 
$$+ \sum_{a=1}^{2}({\kappa}\wedge{\chi_{\boldkey F}})(X_{a},X_{1},D_{X_{a}}X_{2},V_{1},\dots,V_{p-2})
 + \sum_{a=1}^{2}({\kappa}\wedge{\chi_{\boldkey F}})(X_{a},X_{1},X_{2},D_{X_{a}}V_{1},\dots,V_{p-2})$$ 
$$ +\dots+ \sum_{a=1}^{2}({\kappa}\wedge{\chi_{\boldkey F}})(X_{a},X_{1},X_{2},V_{1},\dots,D_{X_{a}}V_{p-2})
%
%
- \sum_{i=1}^{p-2}V_{i}({\kappa}\wedge{\chi_{\boldkey F}})(V_{i},X_{1},X_{2},V_{1},\dots,V_{p-2})$$ 
$$ + \sum_{i=1}^{p-2}({\kappa}\wedge{\chi_{\boldkey F}})(D_{V_{i}}V_{i},X_{1},X_{2},V_{1},\dots,V_{p-2})
+\sum_{i=1}^{p-2}({\kappa}\wedge{\chi_{\boldkey F}})(V_{i},D_{V_{i}}X_{1},X_{2},V_{1},\dots,V_{p-2})$$ 
$$ +\sum_{i=1}^{p-2}({\kappa}\wedge{\chi_{\boldkey F}})(V_{i},X_{1},D_{V_{i}}X_{2},V_{1},\dots,V_{p-2})
+ \sum_{i=1}^{p-2}({\kappa}\wedge{\chi_{\boldkey F}})(V_{i},X_{1},X_{2},D_{V_{i}}V_{1},\dots,V_{p-2})$$ 
$$ + \dots+ Ê\sum_{i=1}^{p-2}({\kappa}\wedge{\chi_{\boldkey F}})(V_{i},X_{1},X_{2},V_{1},\dots,D_{V_{i}}V_{p-2})$$ 
$$ - V_{p-1}({\kappa}\wedge{\chi_{\boldkey F}})(V_{p-1},X_{1},X_{2},V_{1},V_{2},\dots,V_{p-2})
+ ({\kappa}\wedge{\chi_{\boldkey F}})(D_{V_{p-1}}V_{p-1},X_{1},X_{2},V_{1},V_{2},\dots,V_{p-2})$$ 
$$ + ({\kappa}\wedge{\chi_{\boldkey F}})(V_{p-1},D_{V_{p-1}}X_{1},X_{2},V_{1},V_{2},\dots,V_{p-2}) 
+ ({\kappa}\wedge{\chi_{\boldkey F}})(V_{p-1},X_{1},D_{V_{p-1}}X_{2},V_{1},V_{2},\dots,V_{p-2}) $$ 
$$ + ({\kappa}\wedge{\chi_{\boldkey F}})(V_{p-1},X_{1},X_{2},D_{V_{p-1}}V_{1},V_{2},\dots,V_{p-2})
+ \dots + ({\kappa}\wedge{\chi_{\boldkey F}})(V_{p-1},X_{1},X_{2},V_{1},V_{2},\dots, D_{V_{p-1}}V_{p-2})$$ 
$$ - V_{p}({\kappa}\wedge{\chi_{\boldkey F}})(V_{p},X_{1},X_{2},V_{1},V_{2},\dots,V_{p-2})
+ ({\kappa}\wedge{\chi_{\boldkey F}})(D_{V_{p}}V_{p},X_{1},X_{2},V_{1},V_{2},\dots,V_{p-2})$$ 
$$ + ({\kappa}\wedge{\chi_{\boldkey F}})(V_{p},D_{V_{p}}X_{1},X_{2},V_{1},V_{2},\dots,V_{p-2}) 
+ ({\kappa}\wedge{\chi_{\boldkey F}})(V_{p},X_{1},D_{V_{p}}X_{2},V_{1},V_{2},\dots,V_{p-2})  $$ 
$$ + ({\kappa}\wedge{\chi_{\boldkey F}})(V_{p},X_{1},X_{2},D_{V_{p}}V_{1},V_{2},\dots,V_{p-2})
 + \dots + ({\kappa}\wedge{\chi_{\boldkey F}})(V_{p},X_{1},X_{2},V_{1},V_{2},\dots, D_{V_{p}}V_{p-2}).$$ 
\bigskip 
\noindent 
All terms above with two horizontal vector field arguments vanish. Terms in 
\newline\noindent 
 $\sum_{i=1}^{p-2}({\kappa}\wedge{\chi_{\boldkey F}})(V_{i},D_{V_{i}}X_{1},X_{2},V_{1},\dots,V_{p-2})$ 
 vanish individually because the arguments $V_{i}$ repeat when \newline\noindent $1 \le i \le p-2$. Accordingly, 
the only non-zero terms are: 
$$ \ ({\kappa}\wedge{\chi_{\boldkey F}})(V_{p-1},D_{V_{p-1}}X_{1},X_{2},V_{1},V_{2},\dots,V_{p-2})\tag{1.24}$$ 
$$+ ({\kappa}\wedge{\chi_{\boldkey F}})(V_{p-1},X_{1},D_{V_{p-1}}X_{2},V_{1},V_{2},\dots,V_{p-2})$$ 
$$+ ({\kappa}\wedge{\chi_{\boldkey F}})(V_{p},D_{V_{p}}X_{1},X_{2},V_{1},V_{2},\dots,V_{p-2})$$ 
$$+ ({\kappa}\wedge{\chi_{\boldkey F}})(V_{p},X_{1},D_{V_{p}}X_{2},V_{1},V_{2},\dots,V_{p-2}) .$$ 

\bigskip 
\noindent 
Only the vertical components of $D_{V_{l}}X_{j}$ matter in the above calculations because 
when ${\kappa}\wedge{\chi_{\boldkey F}}$ is evaluated on $p+1$ arguments with two or more horizontal the result is zero. 
Recall, ${\Cal V}D_{V_{j}}X_{i} = T_{V_{j}}X_{i}$. Hence, we have, 

$$ ({\kappa}\wedge{\chi_{\boldkey F}})(V_{p-1},T_{V_{p-1}}X_{1},X_{2},V_{1},V_{2},\dots,V_{p-2})\tag{1.25}$$ 
$$+ ({\kappa}\wedge{\chi_{\boldkey F}})(V_{p-1},X_{1},T_{V_{p-1}}X_{2},V_{1},V_{2},\dots,V_{p-2})$$ 
$$+ ({\kappa}\wedge{\chi_{\boldkey F}})(V_{p},T_{V_{p}}X_{1},X_{2},V_{1},V_{2},\dots,V_{p-2})$$ 
$$+ ({\kappa}\wedge{\chi_{\boldkey F}})(V_{p},X_{1},T_{V_{p}}X_{2},V_{1},V_{2},\dots,V_{p-2}) .$$ 
\bigskip 
\noindent 
A routine argument using the properties of the tensor $T$ introduced in the beginning shows 
that \newline\noindent $g(T_{V_{p-1}}X_{j},V_{p}) = g(T_{V_{p}}X_{j},V_{p-1})$, where $j=1$ or $j=2$. This means 
if the $V_{p}$ -component of $T_{V_{p-1}}X_{1}$ is $a$, then the $V_{p-1}$ -component of $T_{V_{p}}X_{1}$ is also $a$. Likewise, 
if $c$ is the $V_{p}$ -component of $T_{V_{p-1}}X_{2}$, then c is also the $V_{p-1}$ -component of $T_{V_{p}}X_{2}$. 
Hence, (1.25) becomes: 

$$ \leqno (1.26)\qquad ({\kappa}\wedge{\chi_{\boldkey F}})(V_{p-1},aV_{p} ,X_{2},V_{1},V_{2},\dots,V_{p-2})$$ 
$$+ ({\kappa}\wedge{\chi_{\boldkey F}})(V_{p-1},X_{1}, cV_{p},V_{1},V_{2},\dots,V_{p-2})$$ 
$$+ ({\kappa}\wedge{\chi_{\boldkey F}})(V_{p},aV_{p-1},X_{2},V_{1},V_{2},\dots,V_{p-2})$$ 
$$+ ({\kappa}\wedge{\chi_{\boldkey F}})(V_{p},X_{1},cV_{p-1},V_{1},V_{2},\dots,V_{p-2}) = 0 .$$ 
\newline 
\noindent 
The proof of Theorem 1.4 is now complete, provided we observe that in the very special case that ${\boldkey F}$ is a flow 
on $M^{3}$, the third computation is superfluous. 
\bigskip 
\bigskip 

\noindent 
As an application of Theorem 1.4 we establish the following result which also 
uses a result of Ranjan (see [Ra]). 
We will follow the exposition of 
Ranjan's Theorem as given in [T3], pages 76 and 77 (see also the Corollary on 
the top of page 89 in [Ra]). 
Essentially, our result says gives necessary and sufficient  conditions for a 
closed $3$-manifold with a Dominguez metric to admit a non-trivial 
local de-Rham decomposition. 
 $Ric^{M}(E,E)$ denotes the Ricci tensor with respect to the Levi-Civita 
connection on $M$ 
evaluated on a vector field $E$.   
\bigskip 
\bigskip

\noindent
\underbar{Corollary 1.5}. Let $M^{3}$ be a closed, oriented, $C^{\infty}$, $3$-manifold, with a transversely oriented Riemannian flow,
$\boldkey F$. Suppose $g$ is a Dominguez metric for the flow , ${\boldkey F}$, and let $V$ be a unit length vector field
tangent to this flow.
\newline
\newline
\noindent If $Ric^{M}(V,V) \equiv 0$ on $M^{3}$ and ${\kappa \wedge \chi_{\boldkey F}}$ 
is harmonic, then ${\boldkey H}$ is integrable,
${\boldkey F}$ is totally geodesic, ${\kappa \wedge \chi_{\boldkey F}} \equiv 0$, 
and locally $M^{3}$ is isometric to a product of the plaques
of the leaves of ${\boldkey H}$ and ${\boldkey F}$.
\newline
\newline
\noindent
Conversely, if ${\boldkey H}$ is integrable and ${\boldkey F}$ is totally geodesic, then on $M^{3}$, $Ric^{M}(V,V) \equiv 0 $ and
${\kappa \wedge \chi_{\boldkey F}} \equiv 0 $. In particular, ${\kappa \wedge \chi_{\boldkey F}}$ is harmonic.
\newline
\newline
\noindent \underbar{Proof}. In the proof we will let $ \lbrace X_{1}, ..., X_{q} \rbrace$ 
denote a local basic orthonormal frame for ${\boldkey H}$,
with $q=2$. We use this seemingly cumbersome notation because the same work will yield 
another somewhat more general result  essentially at
no extra cost. First note, ${\kappa}\wedge \chi_{\boldkey F}$ is closed because of Theorem 1.4 
and the theorem of Kamber-Tondeur ([T3] page 82)
which applies in the case of a Dominguez metric.  Note all the calculations are independent of the 
local orthonormal frame for ${\boldkey H}$.
The idea is to exploit equations 6.22 and 6.21 of [T3] in that order.  Equation 6.22 of [T3] yields
$Ric^{M}(V,V) = div_{M}{\tau} + \sum_{i = 1}^{q} g(A_{X_{i}}V, A_{X_{i}}V)$. If $ Ric^{M}(V,V) = 0 $ 
on $M^{3}$, then integration yields
$ \int_{M}  \sum_{i = 1}^{q} g(A_{X_{i}}V, A_{X_{i}}V) = 0$, so each $A_{X_{i}}V = 0$ so 
$A \equiv 0$. Since $A_{X}Y = (1/2) \Cal{V}[X,Y]$,
we see ${\boldkey H}$ is integrable. Since $A_{X}Y = \Cal{V}D_{X}Y $, ${\boldkey H}$ is totally geodesic.
\newline
\newline
\noindent
Using the fact that for our flow $\Cal{V}D_{X}V = 0$ for any basic $X$, equation 6.21 of [T3] yields \newline
\noindent $Ric^{M}(V,V) = div_{\boldkey H}{\tau} - \sum_{i = 1}^{q}g(T_{V}X_{i}, T_{V}X_{i}) +
 \sum_{i = 1}^{q}g(A_{X_{i}}V, A_{X_{i}}V)$. 
\newline \noindent Now $Ric^{M}(V,V) \equiv 0$ and $A \equiv 0$ and, when ${\kappa}\wedge {\chi}_{\boldkey F}$ is harmonic,
$div_{\boldkey H}{\tau} = 0 $, when $q = 2$. Hence, $\sum_{i = 1}^{q}g(T_{V}X_{i}, T_{V}X_{i}) = 0$. 
This means $T \equiv 0$ or ${\boldkey F}$ is totally geodesic
and so locally $M^{3}$ is isometric to a product of the plaques of the foliations ${\boldkey H}$ and ${\boldkey F}$.
The proof of the converse follows by observing that under the stated hypotheses, 
our version of 6.22 of [T3] yields that $Ric^{M}(V,V) \equiv 0$. Then
our version of 6.21 of [T3] yields $div_{\boldkey H}{\tau}= 0$ which when $q=2$ means, 
${\kappa}\wedge {\chi}_{\boldkey F}$ is harmonic. If $T = 0$,
${\tau} \equiv 0$ and so ${\kappa} \wedge {\chi}_{\boldkey F} \equiv 0$.
\newline
\newline
The second result using the proof above works for a Riemannian flow of arbitrary codimension on a closed, 
connected manifold. 
\newline
\newline
\underbar{Corollary 1.6}. Let $M$ be a closed, oriented, $C^{\infty}$, $n$-manifold, 
with a transversely oriented Riemannian flow,
${\boldkey F}$. Suppose $g$ is a Dominguez metric for the flow , ${\boldkey F}$, and let $V$ be a unit length vector field
tangent to this flow. 
\newline
\newline
\noindent If $Ric^{M}(V,V) \equiv 0$ on $M$ and $div_{\boldkey H}{\tau} = 0$, then ${\boldkey H}$ is integrable,
${\boldkey F}$ is totally geodesic, ${\tau} = 0$, and locally $M$ is isometric to a product of the plaques
of the leaves of ${\boldkey H}$ and ${\boldkey F}$.
\newline
\newline
\noindent
Conversely, if ${\boldkey H}$ is integrable and ${\boldkey F}$ is totally geodesic, 
then on $M$, $Ric^{M}(V,V) \equiv 0 $ and
$div_{\boldkey H} {\tau} = 0 $. In particular, ${\tau} = 0$.
\newline
\newline
\underbar{Proof}. As noted, the result follows from the proof of 1.5 with minor modifications.  
\newline
\newline
\noindent 
Now suppose $M^{n}$ is a closed, connected, oriented, Riemannian manifold admitting a codimension-one Riemannian foliation 
${\boldkey F}$. Let $g$ be a Dominguez metric for ${\boldkey F}$. 
Then $ {\tau} = Ê\sum_{i = 1}^{n-1} {\Cal H}D_{V_{i}}V_{i} $. 
(As above, we follow the conventions in foliations and suppress the 
usual constant.) And $ div_{\boldkey H}{\tau} = g(X, D_{X}{\tau})$, 
where $X$ is a unit length basic vector field. We have 
the following result. 
\bigskip 
\bigskip 
\noindent\underbar{Remark}. For a general codimension-one, 
transversely oriented foliation on a closed, oriented, Riemannian manifold,
Kamber and Tondeur have shown that the leaves of the foliation are minimal submanifolds with respect to the given metric
if and only if $d{\chi}_{\boldkey F} = 0$ as shown in Theorem 7.35 of [T1], page 92. 
But for foliations of codimension one,
it is also the case (see [T1], page 80) that  $d{\chi}_{\boldkey F} = - {\kappa} \wedge {\chi}_{\boldkey F}$. Hence, 
the leaves of the foliation are minimal in this setting if and only if ${\kappa} \wedge {\chi}_{\boldkey F} = 0$, 
or equivalently in this setting, ${\kappa} \wedge {\chi}_{\boldkey F}$ is a harmonic $n$-form, by the Hodge Theorem.
The next result gives a sufficient explicit condition for 
${\kappa} \wedge {\chi}_{\boldkey F}$ to be harmonic in the very special case 
that the codimension-one foliation is bundle-like with respect to a Dominguez metric. 
We include it because the key condition is 
essentially the same as that for the codimension $q=2$ case in Theorem 1.4.
\newline
\newline
\noindent 
\underbar{Theorem 1.7}. Let $M^{n}$ be a closed, connected, $C^{\infty}$, 
oriented Riemannian manifold admitting a transversely oriented, 
codimension-one, Riemannian foliation ${\boldkey F}$. 
Let $g$ be a Dominguez metric for the foliation ${\boldkey F}$ (with respect to 
which Ê${\boldkey F}$ is bundle-like and ${\kappa}$ is basic). 
Then ${\kappa}\wedge{\chi_{\boldkey F}}$ is harmonic (and hence 
by the above remark identically 0 in this case) if and only if  
$div_{\boldkey H}{\tau} = 0$, where ${\tau}$ is the mean 
curvature one-form dual to ${\kappa}$. 
\newline 
\newline 
\newline 
\noindent 
\underbar{Proof}. ${\kappa}\wedge{\chi_{\boldkey F}}$ is an $n$-form and hence is closed. 
Because the chosen metric, $g$, is a Dominguez 
metric, the mean curvature one-form ${\kappa}$ is basic. 
Just as before, it is closed by the Kamber-Tondeur Theorem. 
We will show under the stated hypotheses, $ {\delta}({\kappa}\wedge{\chi_{\boldkey F}}) = 0$. 
We choose a 
{preferred orthonormal frame} at $x \in M$. 
That is, we choosed $\lbrace X, V_{1}, V_{2}, \dots, V_{n-1}\rbrace$, so $X$ is basic, 
with $({\Cal H }D_{X}X)_{x} = 0 $ and so $({\Cal V}D_{V_i}V_{j})_{x} = 0$, 
where $\lbrace V_{1}, V_{2}, \dots, V_{n-1} \rbrace$ is 
an orthonormal frame for ${\boldkey V}$. Then at $ x \in M$, we have, 
$$\alignat{2} 
 & (1.27) \ \ 
\delta({\kappa}\wedge{\chi_{\boldkey F}})&&\negmedspace\negmedspace\negmedspace (V_{1},V_{2}, \dots, V_{n-1}) \\ 
& & = & Ê 
- (D_X({\kappa}\wedge{\chi_{\boldkey F}}))(X,V_{1},V_{2}, \dots, V_{n-1}) 
- \sum_{j =1}^{n-1} (D_{V_j}({\kappa}\wedge{\chi_{\boldkey F}}))(V_{j},V_{1},V_{2}, \dots, V_{n-1})\\ 
 & ÊÊ& =& -X({\kappa}\wedge{\chi_{\boldkey F}})(X,V_{1},V_{2},\dots,V_{n-1}) + 
({\kappa}\wedge{\chi_{\boldkey F}})(D_{X}X,V_{1},V_{2},\dots,V_{n-1})\\ 
 & & & + ({\kappa}\wedge{\chi_{\boldkey F}})(X,D_{X}V_{1},V_{2},\dots,V_{n-1}) 
+\dots+ ({\kappa}\wedge{\chi_{\boldkey F}})(X,V_{1},V_{2}, \dots,D_{X}V_{n-1})\\ 
&&& - \sum_{j=1}^{n-1}V_{j}({\kappa}\wedge{\chi_{\boldkey F}})(V_{j},V_{1},V_{2},\dots,V_{n-1}) 
+ ÊÊ\sum_{j=1}^{n-1}({\kappa}\wedge{\chi_{\boldkey F}})(D_{V_{j}}V_{j},V_{1},V_{2},\dots,V_{n-1})\\ 
&&& + ({\kappa}\wedge{\chi_{\boldkey F}})(V_{1},D_{V_{1}}V_{1},V_{2},\dots,V_{n-1}) + 
({\kappa}\wedge{\chi_{\boldkey F}})(V_{1},V_{1},D_{V{1}}V_{2},\dots,V_{n-1})\\ 
&&& +\dots + ({\kappa}\wedge{\chi_{\boldkey F}})(V_{1},V_{1},V_{2},\dots,D_{V_{1}}V_{n-1})\\ 
&&& + ({\kappa}\wedge{\chi_{\boldkey F}})(V_{2},D_{V{2}}V_{1},V_{2},\dots,V_{n-1}) 
+ ({\kappa}\wedge{\chi_{\boldkey F}})(V_{2},V_{1},D_{V{2}}V_{2},\dots,V_{n-1})\\ 
&&&+ \dots + ({\kappa}\wedge{\chi_{\boldkey F}})(V_{2},V_{1},V_{2},\dots,D_{V_{2}}V_{n-1}) 
+ \dots + ({\kappa}\wedge{\chi_{\boldkey F}})(V_{n-1},D_{V_{n-1}}V_{1},V_{2},\dots,V_{n-1})\\ 
&&&+ \dots + ({\kappa}\wedge{\chi_{\boldkey F}})(V_{n-1},V_{1},V_{2},\dots,D_{V_{n-1}}V_{n-1}).\\ 
\endalignat$$ ÊÊÊÊÊ 
\newline 
\newline 
\noindent 
Just as in the proof of Theorem 1.4, the expressions 
$ \sum_{j=1}^{n-1}({\kappa}\wedge{\chi_{\boldkey F}})(D_{V_{j}}V_{j},V_{1},V_{2},\dots,V_{n-1})$, 
\newline 
\noindent 
$({\kappa}\wedge{\chi_{\boldkey F}})(V_{1},D_{V_{1}}V_{1},V_{2},\dots,V_{n-1})$, 
$({\kappa}\wedge{\chi_{\boldkey F}})(V_{2},V_{1},D_{V{2}}V_{2},\dots,V_{n-1})$, $\dots$, and, 
\newline\noindent 
$({\kappa}\wedge{\chi_{\boldkey F}})(V_{n-1},V_{1},V_{2},\dots,D_{V_{n-1}}V_{n-1})$ sum to zero. 
Except for the first term, the remaining terms in (1.26) vanish because of repeated arguments, the fact 
that $D_{X}V_{j}$ has no non-zero $V_{j}$ component and so these expressions are evaluated with two purely horizontal 
arguments and hence vanish as well. 
\newline 
\newline 
\noindent 
Then (1.27) becomes, 
$$ \,\, -X({\kappa}\wedge{\chi_{\boldkey F}})(X,V_{1},V_{2},\dots,V_{n-1}) 
= -X{\kappa}(X) = -div_{\boldkey H}{\tau} , \tag{1.28}$$ 
which must vanish identically if ${\kappa}\wedge{\chi_{\boldkey F}}$ is co-closed. 
Our theorem will be proven if we can show 
\newline 
\noindent 
$\delta ({\kappa}\wedge{\chi_{\boldkey F}})(X,V_{1},V_{2},\dots, {\hat V_{l}},\dots,V_{n-1}) 
= 0$. Essentially without loss of generality, we will show 
\newline 
\noindent 
$\delta ({\kappa}\wedge{\chi_{\boldkey F}})(X,V_{1},V_{2},\dots,V_{n-2}) 
=0$, since by renumbering the vertical vectors, up to sign, the computation will always evaluate to zero. 
$$\alignat{2} 
 & (1.29) \negmedspace &&\delta 
({\kappa} 
\wedge{\chi_{\boldkey F}}) 
(X,V_{1},V_{2},\dots,V_{n-2})\\ 
&& =& -(D_{X}({\kappa}\wedge{\chi_{\boldkey F}}))(X,X,V_{1},V_{2},\dots,V_{n-2}) 
- \sum_{j=1}^{n-2}(D_{V_{j}}({\kappa}\wedge{\chi_{\boldkey F}}))(V_{j},X,V_{1},V_{2},\dots,V_{n-2})\\ 
&&& -(D_{V_{n-1}}({\kappa}\wedge{\chi_{\boldkey F}}))(V_{n-1},X,V_{1},V_{2},\dots,V_{n-2})\\ 
&& =& - X ({\kappa}\wedge{\chi_{\boldkey F}})(X,X,V_{1},V_{2},\dots,V_{n-2}) 
+ Ê({\kappa}\wedge{\chi_{\boldkey F}})(D_{X}X,X,V_{1},V_{2},\dots,V_{n-2})\\ 
&&& + ({\kappa}\wedge{\chi_{\boldkey F}})(X,D_{X}X,V_{1},V_{2},\dots,V_{n-2}) 
+ ({\kappa}\wedge{\chi_{\boldkey F}})(X,X,D_{X}V_{1},V_{2},\dots,V_{n-2})\\ 
&&& + \dots + ({\kappa}\wedge{\chi_{\boldkey F}})(X,X,V_{1},V_{2},\dots,D_{X}V_{n-2})\\ 
&&& - \sum_{j=1}^{n-2}V_{j}({\kappa}\wedge{\chi_{\boldkey F}})(V_{j},X,V_{1},V_{2},\dots,V_{n-2}) 
 + \sum_{j= 1}^{n-2} ({\kappa}\wedge{\chi_{\boldkey F}})(D_{V_{j}}V_{j},X,V_{1},V_{2},\dots,V_{n-2})\\ 
&&& \negmedspace + \sum_{j=1}^{n-2} ({\kappa}\wedge{\chi_{\boldkey F}})(V_{j},D_{V_{j}}X,V_{1},\dots,V_{j},\dots,V_{n-2}) 
+ \dots + \sum_{j=1}^{n-2}(\sum_{i=1}^{n-2}(\kappa\wedge{\chi_{\boldkey F}})(V_{j},X,\dots,D_{V_{j}}V_{i},\dots,V_{n-2}))\\ 
&&& - V_{n-1}({\kappa}\wedge{\chi_{\boldkey F}})(V_{n-1},X,V_{1},\dots,V_{n-2}) + 
({\kappa}\wedge{\chi_{\boldkey F}})(D_{V_{n-1}}V_{n-1},X,V_{1},\dots,V_{n-2}) \\ 
&&& + ({\kappa}\wedge{\chi_{\boldkey F}})(V_{n-1},D_{V_{n-1}}X,V_{1},\dots,V_{n-2})\\ 
&&&+ ({\kappa}\wedge{\chi_{\boldkey F}})(V_{n-1},X,D_{V_{n-1}}V_{1},\dots,V_{n-2}) 
 + \dots + ({\kappa}\wedge{\chi_{\boldkey F}})(V_{n-1},X,V_{1},\dots,D_{V_{n-1}}V_{n-2}).\\ 
\endalignat $$ 
\noindent Then (1.29) becomes, 
$$ -V_{n-1}({\kappa}\wedge{\chi_{\boldkey F}})(V_{n-1},X,V_{1},\dots,V_{n-2}) = \pm V_{n-1}{\kappa}(X) 
 = \pm d{\kappa}(V_{n-1},X) = 0 ,  \tag{1.30}$$ 
because ${\kappa}$ is a basic form when $g$ is a Dominguez metric. This completes the proof of Theorem 1.7. 
\newline 
\newline 
\noindent 
Following [T3], page 99, we define the following connection, ${\tilde D}$. For vector fields $E$ and $F$ on $M$, we set:

$${\tilde D}_{E}F = {\Cal V}D_{E}{\Cal V}F+ {\Cal H}D_{E}{\Cal H}F , \tag{1.31} $$ 

\noindent 
where $D$ is the Levi-Civita connection on $M$. ÊAgain following [T3] (page 102) or [Mi-Ri-To], 
let ${\omega}$ be a basic $r$-form. Let 
$ \lbrace E_{2},\dots,E_{r} \rbrace$ be vector fields on $M$. Let $\lbrace V_{1},\dots,V_{p},X_{1},\dots,X_{q} \rbrace$ 
be an orthonormal frame for a bundle-like foliation of leaf dimension p and codimension q. Set, 

$$ 
\alignat{3}
&  {\tilde {\delta}}{\omega} (E_{2},\dots,E_{r}) & = &
- \sum_{j = 1}^{p} V_{j}({\omega}(V_{j},E_{2},\dots,E_{r}) 
 + \sum_{j = 1}^{p} {\omega}({\tilde {D}}_{V_{j}}V_{j}, E_{2},\dots,E_{r}) \tag{1.32}\\ 
& & & + \sum_{j=1}^{p} \sum_{i=2}^{r}{\omega}(V_{j},E_{2},\dots, {\tilde {D}}_{V_{j}}E_{i},\dots,E_{r}) 
- \sum_{k=1}^{q}X_{k}({\omega}(X_{k},E_{2},\dots,E_{r}))\\
& & & + \sum_{k=1}^{q}{\omega}({\tilde {D}}_{X_{k}}X_{k},E_{2},\dots,E_{r}) + 
\sum_{k=1}^{q}\sum_{i=2}^{r}{\omega}(X_{k},E_{2},\dots, {\tilde {D}}_{X_{k}}E_{i}, \dots, E_{r}).
\endalignat 
$$ 

\noindent 
Then by [T3], page 102, if ${\omega}$ is a basic $r$-form, ${\tilde {\delta}}{\omega}$ is a basic $(r-1)$-form. In particular, 
if ${\boldkey F}$ is a transversely oriented, Riemannian foliation on a closed, oriented Riemannian manifold $(M,g)$, where 
$g$ is a Dominguez metric for ${\boldkey F}$, a straightforward calculation yields the following: 

$$ {\tilde {\delta}}{\kappa} = - div_{\boldkey {H}}{\tau} . \tag{1.33}$$ 
\noindent 
We have the following theorem which combines Theorems 1.4 and 1.7. 
\newline 
\newline 
\noindent 
\underbar{Theorem 1.8}. Let $(M^{n},g) $ be a closed, oriented, $C^{\infty}$, Riemannian manifold, 
with a transversely oriented, 
codimension-q, Riemannian foliation ${\boldkey F}$, with $q=1$ or $q=2$. 
Suppose $g$ is a Dominguez metric for ${\boldkey F}$. 
Then ${\kappa}\wedge{\chi_{\boldkey F}}$ is harmonic on $M$ if and only if 
${\tilde {{\delta}}{\kappa} }= 0$. Under the stated hypotheses
when $q=1$, ${\kappa}\wedge{\chi_{\boldkey F}} = 0$.
\newline 
\newline 
\noindent 
\underbar{Proof}. The proof follows immediately from Theorems 1.4 and 1.7 and the above remarks. 
\newline 
\newline 
\noindent 
It might be useful to rephrase Theorem 1.8 in the following way. 
It should be noted however, that thanks to the fundamental result of 
Dominguez, we always know there exists a bundle-like metric $g$ for ${\boldkey F}$ so ${\kappa}$ is basic, 
so in a sense the reformulation is redundant. 
\newline 
\newline 
\noindent 
\underbar{Theorem 1.9}. Let $(M^{n},g) $ be a closed, oriented, 
$C^{\infty}$, Riemannian manifold, with a transversely oriented, 
codimension-q, foliation ${\boldkey F}$, with $q=1$ or $q=2$. Suppose ${\boldkey F}$ is bundle-like with respect to Ê$g$. 
Then ${\kappa}\wedge{\chi_{\boldkey F}}$ is harmonic on $M$ if and only if ${\kappa}$ is basic and 
${\tilde {{\delta}}{\kappa}} = 0$. Under the stated hypotheses when $q=1$, ${\kappa}\wedge{\chi_{\boldkey F}} = 0$. 
\newline 
\newline 
\noindent 
\underbar{Proof}. ${\kappa}\wedge{\chi_{\boldkey F}}$ is always closed if $q=1$. 
If $q=2$, then ${\kappa}\wedge{\chi_{\boldkey F}}$ 
is closed if and only if Ê${\kappa}$ is horizontally closed by (1.14). But if ${\kappa}$ is basic, ${\kappa}$ is closed 
by the already mentioned result of Kamber-Tondeur. ${\delta}({\kappa}\wedge{\chi_{\boldkey F}}) = 0$ iff Ê 
$\tilde{\delta}{\kappa} = - div_{\boldkey H}{\tau}= 0$. 
\newline 
\newline 
\noindent 
\underbar{Remark}. A straightforward calculation shows ${\delta}{\kappa} = {\kappa}({\tau}) - div_{\boldkey{H}}{\tau}$. 
If additionally, Ê${\delta}{\kappa} = 0$, then ${\kappa}$ would be closed and co-closed and hence harmonic 
on $M$ itself, a situation 
not necessary to our work here. 
If we set ${\tilde {\triangle}} = d{\tilde {\delta}} + {\tilde {\delta}}d$ as in [T3], page 102, 
then $d{\kappa} = 0$ and ${\tilde{\delta}}{\kappa} = 0$ implies ${\tilde {\triangle}}{\kappa} = 0$. 
However, ${\tilde {\triangle}}$ is not self-adjoint. 
\newline
\newline

\noindent Now let $M^{n}$ be any closed, oriented Riemannian manifold admitting 
a transversely oriented foliation ${\boldkey F}$ of leaf dimension p and 
codimension q. Let ${\tau}$ be the mean curvature vector field of the foliation ${\boldkey F}$.

\noindent Then,

$$ div_{M}{\tau}  = \Sigma_{ {\alpha}  = 1}^{n-p} g(D_{X_{\alpha}} {\tau}, X_{\alpha}) 
+ \Sigma_{i =1}^{p} g(D_{V_{i}}{\tau}, V_{i}).\tag{1.34}$$
\newline
\newline
\noindent Using the standard properties of the tensor $T$, this becomes

$$div_{M}{\tau} + g({\tau},{\tau}) = div_{\boldkey H}{\tau}.\tag{1.35} $$
\newline
\newline
\noindent Integrating we get,

$$\int_{M} g({\tau},{\tau}) dV =  \int_{M} div_{\boldkey H} {\tau} dV.\tag{1.36} $$
\newline
\newline

\noindent We have the following lemma, which applies to an arbitrary transversely oriented foliation 
on a closed oriented Riemannian manifold $M^{n}$, 
not just Riemannian foliations.
\newline
\newline

\noindent \underbar{Lemma 1.10}. Let ${\boldkey F}$ be any transversely oriented foliation of leaf dimension p on a closed, 
oriented, Riemannian manifold
$M^{n}$. 
\newline
\newline
\noindent (1) If $ \int_{M} div_{\boldkey H} {\tau} dV = 0$, then ${\tau} = 0$ and the leaves of ${\boldkey F}$ are minimal.
\newline
\newline
\noindent (2) Conversely, if ${\tau} = 0$ on such an $M^{n}$, 
then  $div_{\boldkey H} {\tau} = 0$ so $\int_{M} div_{\boldkey H} {\tau} dV = 0  .$
\newline
\newline
\noindent \underbar{Theorem 1.11}. Let $M^{n}$ be a closed, 
oriented Riemannian manifold admitting a transversely oriented Riemannian foliation 
${\boldkey F}$ of codimension q with $q=1$ or $2$. Then ${\kappa} \wedge {\chi}_{\boldkey F}$ 
is harmonic with respect to a Dominguez metric for ${\boldkey F}$
if and only if the mean curvature one-form ${\kappa} = 0$ and so the 
leaves of ${\boldkey F}$ are minimal submanifolds of $M^{n}$.
\newline
\newline
\noindent \underbar{Proof}. For $q = 2$ or $1$,  Theorems 1.4 and 1.7 respectively  
guarantee ${\kappa} \wedge {\chi}_{\boldkey F}$ is harmonic if and only if 
$div_{\boldkey H} {\tau} = 0$. The result now follows directly from Lemma 1.10. 
 
\vskip .75in 
\noindent {\bf 2. Bundle-like foliations with totally umbilical leaves} 

\vskip .5in 
\noindent 
We begin by recalling some basic local properties of Riemannian submersions 
and of bundle-like foliations. 
The convention for the Riemannian tensor on a Riemannian manifold $(M,g)$ is: 
$$R(E,F)G={D_E}{D_F}G-{D_F}{D_E}G-{D_{[E,F]}}G 
\ \ \text{and} \ \ R(E,F,G,G')=-g(R(E,F)G,G').$$ 

\noindent 
If $F$ is a $p$-dimensional leaf of foliation ${\boldkey F}$ then $T_UV$ is the second fundamental form 
of the leaf and the mean curvature vector field $\tau$ is given by: 
$$\tau=\sum_{i=1}^{p}T_{V_i}V_i,$$ 
where $\{V_i\}_{1\leq i\leq p}$ is a local orthonormal frame of vector fields tangent to leaves. 

\noindent 
A $p$-dimensional submanifold $F$ of a Riemannian manifold $(M,g)$ 
is said to be totally umbilical if the second 
fundamental form $T$ is given by, $T(U,V)=(1/p)g(U,V)\tau$ for any Êvectors $U$, $V$ tangent to ${\boldkey F}$.

\bigskip 
\noindent 
The following equations, usually called {\it O'Neill's equations}, 
characterize the geometry of a bundle-like foliation ${\boldkey F}$ on $(M,g)$ 
(see [T3] page 51, or the known results for Riemannian submersion [O'N, Gr]). 

\bigskip 
\noindent 
\underbar{Proposition 2.1}. 
For every vertical vector fields $U$, $V$, $W$, $W'$ 
and for every horizontal vector fields $X$, $Y$, $Z$, $Z'$, 
we have the following formulas: 
\smallskip 
\item{i)} $R(U,V,W,W')=\hat R(U,V,W,W')-g(T_UW,T_VW')+g(T_VW,T_UW')$, 
\item{ii)} $R(U,V,W,X)=g((D _VT)_UW,X)-g((D _UT)_VW,X),$ 
\item{iii)} $R(X,U,Y,V)=g((D _XT)_UV,Y)-g(T_UX,T_VY)+g((D_UA)_XY,V)+g(A_XU,A_YV)$, 
\item{iv)} Ê$R(X,Y,Z,U)=g((D_ZA)_XY,U)+g(A_XY,T_UZ)-g(A_YZ,T_UX)-g(A_ZX,T_UY)$, 
\item{v)} $R(X,Y,Z,Z')=R^*( X, Y, Z, Z') 
 -2g(A_XY,A_ZZ')+g(A_YZ,A_XZ')-g(A_XZ,A_YZ'),$ 
\smallskip 

\noindent 
where we denote by $R$, $\hat{R}$ and $R^*$ the Riemannian tensors for the connections 
$D$ of $M$, $\hat{D}$ of ${\boldkey F}$, and $D^*$ on the transversal distribution ${\boldkey H}$, 
respectively.

\bigskip 
\noindent 
Using O'Neill's equations, we get the following lemma. 

\bigskip 
\noindent 
\underbar{Lemma Ê2.2}. 
If Ê$\boldkey F$ is a bundle-like foliation on $(M,g)$ 
with totally umbilical leaves then: 
$$\alignat 2 
&a)\ R(U,V,U,V) & = &\hat{R}(U,V,U,V) +[g(U,V)^2-g(U,U)g(V,V)]g(\frac{\tau}{p},\frac{\tau}{p});\\ 
&b)\ R(X,U,X,U) & = &g(U,U)[g(D_X\frac{\tau}{p},X)-g(X,\frac{\tau}{p})^2]+g(A_XU,A_XU);\\ 
&c)\ R(X,Y,X,Y) & = & R^*( X, Y, X,Y)-3g(A_XY,A_XY).\\ 
\endalignat $$ 

\bigskip 
\noindent 
\underbar{Proposition 2.3}. 
Let $(M,g)$ be a Riemannian manifold with a bundle-like foliation ${\boldkey F}$. 
We assume that ${\boldkey F}$ has totally umbilical leaves and 
$X$, $Y$ are basic vector fields. 
Then $A_XY$ is a Killing vector field 
 along leaves if and only if $g(D_{X}\tau,Y)=g(D_{Y}\tau,X)$. 
\newline 
\newline 
\noindent 
\underbar{Proof}. 
Using Proposition 2.1 from [EP] we have: 
$$ g(D_U(A_{X}Y),V)+g(D_V(A_{X}Y),U)=g(U,V)d\kappa(X,Y). \tag{2.1}$$ 
On the other hand, 
$$ 
\alignat 2 
& d\kappa(X,Y) & = & Xg(\tau,Y)-g(\tau,D_{X}Y)-Yg(\tau,X)+g(\tau,D_{Y}X)\\ 
&              & = & g(D_{X}\tau,Y)-g(D_{Y}\tau,X).\\ 
\endalignat 
$$

\bigskip 
\noindent 
\underbar{Remark}. (i) The affirmation of Proposition 2.3 holds if $\tau$ is parallel 
in the transversal distribution along leaves. 

\noindent 
(ii) $\kappa$ is horizontally closed if and only if $g(D_{X}\tau,Y)=g(D_{Y}\tau,X)$ 
for any horizontal vector fields $X$, $Y$. 

\bigskip 

\noindent
In the next proposition we establish under some certain conditions that  $\kappa$ is a basic one-form.

\bigskip 

\noindent 
\underbar{Proposition 2.4.} 
Let $(M,g)$ be a Riemannian manifold with constant curvature. If ${\boldkey F}$ is a Ê 
bundle-like foliation with totally umbilical leaves and the dimension of the leaves 
is $p\geq 2$ then: 
\item{a)} $\tau$ is parallel in the 
 transversal distribution 
along leaves. 
(i.e $\Cal{H} D_V\tau=0$ for any vertical vector field $V$). 
\item{b)} $\tau$ is basic, which implies that $\kappa$ is basic.
\item{c)} $A_\tau=0$. 

\medskip 
\noindent \underbar{Proof}. 
Let $\{V_i\}_{1\leq i\leq p}$ be Êa local orthonormal basis of vertical vector fields. 
Let $V$ be Êa vertical vector field and $X$ an horizontal one. 
By Proposition 2.1, we have: 
$$ 
   \sum_{i=1}^{p} ÊÊR(V_i,V,V_i,X)=\sum_{i=1}^{p} g((D_VT)_{V_i}V_i,X)- \sum_{i=1}^{p} g((D_{V_i}T)_V{V_i},X). 
\tag{$*$}
$$ 
Since $(M,g)$ is of constant curvature Êwe have: $\sum_{i=1}^{p} ÊÊR(V_i,V,V_i,X)=0$. 
We compute the first term of the right hand side, 
$$\alignat 2 
&\sum_{i=1}^{p} g((D_VT)_{V_i}V_i,X) & = & \sum_{i=1}^{p} g(D_V(T_{V_i}V_i),X) 
                                         -\sum_{i=1}^{p} g(T_{D_{V}V_i}V_i),X) 
                                         -\sum_{i=1}^{p} g(T_{V_i}(D_V{V_i}),X)\\ 
&                                       & = & Êg(D_{V}\tau ,X) 
                                       -2(1/p) \sum_{i=1}^{p} g(D_{V}V_i,V_i)g(\tau,X)\\ 
&                                     & = & g(D_{V}\tau,X)-(1/p) \sum_{i=1}^{p}Vg(V_i,V_i)g(\tau,X)\\ 
&                                     & = & g(D_{V}\tau,X).\\ 
\endalignat 
$$ 
Then for the second term of the right hand side of $(*)$ we get: 
$$\alignat 2 
&\sum_{i=1}^{p} g((D_{V_i}T)_{V}V_i,X) & = & \sum_{i=1}^{p} g(D_{V_i}T_{V}V_i,X) 
                                          -\sum_{i=1}^{p} g(T_{D_{V_i}V}V_i,X) 
                                         -\sum_{i=1}^{p} g(T_VD_{V_i}V_i,X)\\ 
&                                        & =& (1/p)[ \sum_{i=1}^{p} g(D_{V_i}(g(V,V_i)\tau),X) 
                                     -g(\tau,X)g(D_{V_i}V,V_i)-g(\tau,X)g(V,D_{V_i}V_i)]\\ 
&                                      & =& (1/p)[ \sum_{i=1}^{p} V_i(g(V,V_i)g(\tau,X))-g(V,V_i)g(\tau,D_{V_i}X)) 
                                     -g(\tau,X)V_i(g(V,V_i))]\\ 
&                                      & =& (1/p)g(D_{V}\tau,X).\\ 
\endalignat 
$$ 
Therefore, 
$$(1-(1/p))g(D_{V}\tau,X)=0$$ 
for every vertical vector field $V$ and for every horizontal vector field $X$, 
which Êimplies that $\tau$ is parallel in the transversal distribution 
along leaves. 

\medskip 
\noindent 
b) It is sufficient to show that, for any basic vector field $X$, 
$g({\tau},X)$ is constant along leaves (i.e. $Vg({\tau},X)=0$ for any vertical vector field $V$). 

\noindent 
First we shall establish that $g(A_XY,A_XY)$ is constant along leaves for any basic vector fields $X$, $Y$ . 
Since ${\boldkey F}$ is a bundle-like foliation we can consider a local model $B$ in distinguished 
chart $U$ on $M$. Then the restriction of the foliation ${\boldkey F}$ to $U$ gives a Riemannian submersion 
$\pi:{\boldkey F}/U\to B$ and we have: 
$$R(X,Y,X,Y)=R^*({\pi} X,{\pi}Y,{\pi}X,{\pi}Y)-3g(A_XY,A_XY).$$ 
Since $(M,g)$ is of constant curvature we get 
$g(A_XY,A_XY)$ is constant along leaves for any basic vector fields $X$, $Y$. 
 By polarization, it follows 
$g(A_XY,A_XZ)$ is constant along leaves for any basic vector fields $X$, $Y$, $Z$. 
Therefore, $A_XA_XZ$ is a basic vector field and then again by polarization, 
$A_XA_YZ+A_YA_XZ$ is a basic vector field Êfor any basic vector fields $X$, $Y$, $Z$. 

\smallskip 
\noindent 
By O'Neill equation iv) in Proposition 2.1 we get: 
$$ 
\alignat 2 
& ÊÊÊÊÊR(X,Y,X,A_XY) & = &g((D_XA)_XY,A_XY)+2g(A_XY,T_{A_XY}X)\\ 
& ÊÊÊÊÊÊÊÊÊÊÊÊÊÊÊÊÊÊÊÊ& = &g(D_X A_XY,A_XY)-g(A_{D_XX}Y,A_XY)\\ 
& ÊÊÊÊÊÊÊÊÊÊÊÊÊÊÊÊÊÊÊÊÊÊ& & -g(A_XD_XY,A_XY) +2g(A_XY,T_{A_XY}X)\\ 
  \qquad \qquad & ÊÊÊÊ\qquad \qquad ÊÊÊÊÊÊÊÊÊÊÊÊÊÊ& = &(1/2) X(g(A_XY,A_XY))-g(A_Y{\Cal{H}D_XX},A_YX) 
-g(A_X\Cal{H}D_XY,A_XY) \tag{2.2}\\ 
& ÊÊÊÊÊÊÊÊÊÊÊÊÊÊÊÊÊÊÊÊÊÊÊ& & +2g(A_XY,T_{A_XY}X).\\ 
\endalignat 
$$ 

\noindent 
Since $M$ has constant curvature and $A_XY$ is a vertical vector, we see that 
$ R(X,Y,X,A_XY)=0$. The first three terms of (2.2) are constant along leaves and 
so should be Êthe last one, that means 
$V(g(A_XY,T_{A_XY}X))=0$. 
Since the leaves are totally umbilical, it follows: 
$$\leqno (2.3) \qquad\qquad V(g(X,\tau))g(A_XY,A_XY)=0.$$ 

\noindent 
Let $x\in M$, let $X$, $Y$ be basic vector fields such that 
$Y(x)=\tau(x)$. 
If Ê$(A_X\tau)(x)\not=0$ then, by (2.3), $ V(g(X,\tau))=0$. 
Now we shall consider the case when $(A_X\tau)(x)=0$. 
Then: 
$$ 
 Vg(\tau,X)_x=(g(D_{V}{\tau},X)+g({\tau},D_{V}X))_x= 
g(\tau,A_XV)_x=-g(A_X{\tau},V)_x=0.$$ 
In the second equality we have used a), $D_{V}{\tau}=0$. 

\noindent 
Therefore, $\tau$ is a basic vector field. 

\medskip 
\noindent 
c) The basicity of $\tau$ implies $[V,\tau]$ is a vertical vector field, 
so $0=\Cal{H} [V,\tau]=\Cal{H} D_{V}\tau-\Cal{H} D_{\tau}V$. 
Therefore, 
$$g(A_{\tau}X,V)=g(D_{\tau}X,V)=-g(X,D_{\tau}V)=-g(X,D_{V}{\tau})=0,$$ 
since $\Cal{H} D_{V}\tau=0$. It follows, $A_{\tau}X=0$ for any horizontal vector $X$. 
For a vertical vector $V$, we have, $g(A_{\tau} V,X)=-g(V,A_{\tau}X)=0$ which implies $A_{\tau}V=0$. 
%

\bigskip 
\noindent 
Using part (b) of Lemma 2.2  we get the following proposition.

\bigskip 
\noindent 
\underbar{Proposition 2.5}. 
Let $(M,g)$ be an $n$-dimensional Riemannian manifold with constant curvature $c$. 
Let ${\boldkey F}$ be a bundle-like foliation
with totally umbilical leaves of dimension $p=n-q$ on $(M,g)$.
Then
$$div_{{\boldkey H}}\tau= cpq +\frac{1}{p}g(\tau,\tau)-
p\frac{1}{g(U,U)}\sum_{a=1}^q g(A_{X_a}U,A_{X_a}U), \tag{2.4}$$
for any non-zero vertical vector $U$ and for any orthonormal frame $\{X_a\}$ of
${\boldkey H}$.
\newline 
\newline 
\noindent 
\underbar{Proof}.
By Lemma 2.2 (b) we get
$$
\alignat 2 
& R(X_a,U,X_a,U) &=& g(U,U)(g(D_{X_a}\frac{\tau}p,X_a)-g(X_a,\frac{\tau}p)g(X_a,\frac{\tau}p))
                      +g(A_{X_a}U,A_{X_a}U)\\
& &=& \frac{1}{p} g(U,U)(g(D_{X_a}\tau,X_a)-\frac{1}{p}g(\tau,\tau))+g(A_{X_a}U,A_{X_a}U).\\
\endalignat 
$$ 
\noindent
But, since $M$ has constant curvature $c$, we have
$$
\sum_{a=1}^q  R(X_a,U,X_a,U)=qcg(U,U),
$$
which implies (2.4).

\bigskip 
\noindent
Now if we assume $A\equiv 0$ in Proposition 2.5 we obtain the following result.

\bigskip 
\noindent 
\underbar{Corollary 2.6}. 
Let $(M,g)$ be an $n$-dimensional Riemannian manifold with constant curvature $c$,
and let ${\boldkey F}$ be a bundle-like foliation
with totally umbilical leaves and with horizontal integrable distribution
of dimension $p=n-q$ on $(M,g)$.
If $div_{\boldkey H}\tau= 0$ then
\newline
(i)\ $c\leq 0$ and
\newline
(ii)\ $g(\tau,\tau)=-pqc$.  

\bigskip 
\noindent
Since Proposition 2.4 ensures us that $\kappa$ is a basic one-form,
as a consequence of Theorems 1.4 and 1.7, we get the following
result from Corollary 2.6.

\bigskip 
\noindent 
\underbar{Theorem 2.7}.
Let $(M,g)$ be an oriented closed $n$-dimensional Riemannian manifold with constant curvature $c$,
with $g$ a bundle-like metric for a transversally oriented  foliation $F$ of codimension $q$,
with totally umbilical leaves of dimension $p$.
If the horizontal distribution is integrable and $p>1$ and $q\in\{1,2\}$ then
$\kappa\wedge\chi_{\boldkey F}$ is harmonic if and only if 
$g(\tau,\tau)=-pqc$.

\bigskip 
\noindent 
In Theorem 2.8, we shall consider the case when the horizontal distribution ${\boldkey H}$ is non-integrable
(i.e. $A\not=0$ at every point). 

\bigskip 
\noindent 
\underbar{Theorem 2.8}. 
Let $(M,g)$ be a Riemannian manifold of constant curvature $c$ and let ${\boldkey F}$ be a
bundle-like foliation with totally umbilical leaves on $M$ with 
the leaf dimension 
$p\geq 2$. 
If the transversal distribution ${\boldkey H}$ is non-integrable 
and the transversal curvature operator $R^*$ is 
Einstein-like, then ${\boldkey F}$ is totally geodesic. 

\medskip 
\noindent \underbar{Proof}. 
From O'Neill equations we obtain: 
$$ \sum_{a=1}^q R^*(Y,X_a,Y,X_a)=(q-1)cg(Y,Y) + 3\sum_{a=1}^q g(A_{Y}X_a,A_{Y}X_a), \tag{2.5}$$ 
where $\{X_a\}_{1\leq a\leq q}$ is a local orthonormal frame of ${\boldkey H}={\boldkey V}^\perp$. 
By hypothesis, there exists a basic function $\lambda$ on $M$ such that: 
$$ Ric^*(Y,Y)=\sum_{a=1}^q R^*(Y,X_a,Y,X_a)=\lambda g(Y,Y) \tag{2.6}$$ 
and $\lambda\not=(q-1)c$ since ${\boldkey H}$ is a non-integrable distribution, i.e. 
$\sum_{a=1}^q g(A_{Y}X_a,A_{Y}X_a)>0$ for some $Y$. 

\noindent 
Now taking $Y=\tau$ in (2.5), from Proposition 2.4, $A_{\tau}=0$, and from 
 (2.5) and (2.6), we get: 
$$(\lambda -(q-1)c)g(\tau,\tau)=0,$$ which implies $\tau=0$. 

\bigskip 

\bigskip 

\bigskip 
\noindent 
\underbar{References} 

\bigskip 

\item {[AMR]} R. Abraham, J. E. Marsden, T. Ratiu, 
\underbar{Manifolds, Tensor Analysis and Applications, (Second 
edition)}, \-Sp\-rin\-ger-Verlag, New York, Berlin, Heidelberg, 1988. 

\item {[Bl]} D. Blair, \underbar {Riemannian Geometry of Contact and 
Symplectic Manifolds}, Birkh\"auser, Boston, 2002. 

\item{[B]} A. Borel, {\it Compact Clifford-Klein forms of symmetric 
spaces}, Topology, {\bf 2} (1963), 111--122. 

\item{[C]} G. Cairns, These, {\it Feuilletages geodesibles}, l'Universite 
des Sciences du Languedoc, (1987). 

\item{[CE1]} G. Cairns and R. Escobales, {\it Further geometry of the mean 
curvature one-form and the normal plane field one-form on a foliated 
Riemannian manifold}, J. Australian Math. Soc. (Series A) {\bf 62} (1997), 
46-63. 

\item{[CE2]} G. Cairns and R. Escobales, {\it Note on a theorem of 
Gromoll-Grove}, Bull. 
Austral. Math Soc. {\bf 55} (1997), 1--5. 

\item {[CH]} S.S. Chern and R.S. Hamilton, {\it On Riemannian metrics 
adapted to three-dimensional contact manifolds with Appendix by A. 
Weinstein}, \underbar {Lecture Notes in Mathematics}, {\bf 1111} (Workshop 
Bonn 1984), Springer-Verlag, NY-Berlin (1985), 105--134. 

\item{[D]} D. Dominguez, {\it Finiteness and tenseness theorems 
for Riemannian foliations}, 
Amer. J. Math., {\bf 120} (1998), no. 6, 1237--1276. 

\item{[Gr]} A. Gray, 
{\it Pseudo-Riemannian almost product manifolds and submersions},
J. Math. Mech {\bf 16} (1967), 715--737.

\item{[E1]} R. Escobales, Jr., {\it The integrability tensor for 
bundle-like foliations}, 
Transactions of the American Mathematical Society 
{\bf 270} (1982), 333--339. 

\item{[E2]} R. Escobales, Jr., {\it Foliations by minimal surfaces and 
contact structures on certain closed 3-manifolds}, International 
Journal of Mathematics and Mathematical Sciences, Volume 2003, No.21, 
11 April 2003, 1323--1330. (MR1990563 (2004b:57040)). 

\item{[EP]} R. Escobales and P. Parker, {\it Geometric consequences of the 
normal 
curvature cohomology class in umbilic foliations}, Indiana University 
Mathematics Journal, {\bf 37} (1988), 389--408. 

\item{[JW]} D.L. Johnson and L.B. Whitt,{\it Totally geodesic 
foliations}, Journal of Differential Geometry {\bf 15, 2}, (1980), 225--235. 

\item{[L]} H. B. Lawson, Jr., {\it The quantitative theory of foliations}, 
CBMS 
Regional Conferences in Mathematics, American Mathematical Society 
{\bf 27} (1977).

\item{[Mi-R-To]} M. Min-Oo, E. Ruh, and Ph. Tondeur,{\it Vanishing theorems for the basic cohomology 
of Riemannian foliations}, J. Reine Angew. Math. {\bf 415} (1991), 167--174. 

\item{[O'N]} B. O'Neill, {\it ÊThe fundamental equations of a 
submersion}, Michigan Math. J. {\bf 13} (1966), 459--469. 

\item{[P]} W. Poor, \underbar{Differential Geometric Structures}, 
McGraw-Hill Book Company, New York, 1981. 

\item{[Ra]} A. Ranjan, 
{\it Structural equations and an integral formula for foliated manifolds},
Geom. Dedicata {\bf 20} (1986), 85--91.

\item{[Ru]} H. Rummler, 
 {\it Quelques notions simples en g\'eom\'etrie riemannienne et leurs applications aux 
 Êfeuilletages compacts.} Comment. Math. Helv. {\bf 54} (1979), no. 2, 224--239. 

\item {[S]} P. Schweitzer, S.J., {\it Existence of Codimension One 
Foliations with Minimal Leaves}, Ann. Global Anal. Geom. {\bf 9}, 
{\#1}, (1991), 77--81. 

\item {[Su]} D. Sullivan, {\it A homological characterization of 
foliations consisting of minimal surfaces}, Comment. Math. Helv. 
{\bf 54} (1979), 218--223. 

\item{[T1]} Ph. Tondeur, \underbar{Foliations on Riemannian manifolds}, 
Springer-Verlag, New York, 1988. 

\item{[T2]} Ph. Tondeur, \underbar{Geometry of Riemannian 
Foliations}, Seminar on Mathematical Sciences, No. 20, Keio 
University, Yokohama, Japan, 1994. 

\item{[T3]} Ph. Tondeur, \underbar{Geometry of Foliations}, Birkh\"auser, 
Basel, Boston, Berlin, 1997. 

\item{[Wa]} G. Walschap, {\it Umbilic foliations and curvature}, 
Illinois J. Math. {\bf 41} (1997), no. 1, 122--128.

\vskip .75in 
\noindent 
G. Baditoiu 

\noindent 
Institute of Mathematics of the Romanian Academy and 

\noindent 
Boston University, Department of Mathematics, 111 Cummington St Rm 142, Boston, MA 02215, USA

\noindent 
baditoiu\@math.bu.edu

\vskip .2in 
\noindent 
R.H. Escobales, Jr. 

\noindent 
Canisius College 

\noindent 
Buffalo NY 14208 

\noindent 
escobalr\@canisius.edu

\vskip .2in 
\noindent 
S. Ianus 

\noindent 
University of Bucharest, Department of Mathematics 

\noindent 
P.O. Box 10--119, Bucharest, Romania 

\noindent 
ianus\@gta.math.unibuc.ro

\end